\documentclass{amsart}

\newtheorem{theorem}{Theorem}[section]
\newtheorem{lemma}[theorem]{Lemma}

\theoremstyle{definition}
\newtheorem{definition}[theorem]{Definition}

\theoremstyle{remark}
\newtheorem{remark}[theorem]{Remark}

\numberwithin{equation}{section}

\begin{document}

\title[On the solution of the Collatz problem]{On the solution of the Collatz problem}

\author{Shan-Guang TAN}
\address{
Zhejiang University, Hangzhou, 310027, CHINA}
\curraddr{
Zhejiang University, Hangzhou, 310027, CHINA}
\email{tanshanguang@163.com}



\date{November 17, 2025 and, in revised form, December 29, 2025.}



\begin{abstract}
In this paper, we first prove that given a nonnegative integer $m$ and an odd number $t$ not divisible by $3$, there exists a unique Collatz's Sequence
\[
S_{c}(m,t)=\{n_{0}(m,t),n_{1}(m,t),n_{2}(m,t),\ldots,n_{m}(m,t),n_{m+1}(m,t)\}
\]
produced by a function $n_{i+1}(m,t)=(3n_{i}(m,t)+1)/2$ for $i=0,1,2,\ldots,m$ and ended by an even number $n_{m+1}(m,t)$ where $n_{i}(m,t)=2^{m+1-i}\times3^{i}t-1$ for $i=0,1,2,\ldots,m+1$, by which all odd numbers can be expressed.

Then we discuss the Collatz problem in two ways and prove that each Collatz's Sequence always returns to 1, i.e., the Collatz problem is solved.
\end{abstract}

\maketitle

\section*{Introduction}

A problem posed by L. Collatz in 1937, also called the $3n+1$ problem. Let $m_{0}^{*}$ be an integer. Then one form of Collatz problem asks if iterating
\begin{equation}\label{def:collatz_orig}
m_{k}^{*}=\{\begin{array}{cc}
         m_{k-1}^{*}/2 & \texttt{ for }m_{k-1}^{*}\texttt{ even } \\
        3m_{k-1}^{*}+1 & \texttt{ for }m_{k-1}^{*}\texttt{ odd }
      \end{array}
\end{equation}
always returns to 1 for any positive $m_{0}^{*}$~\cite{N,S}.

The members of the sequence produced by the Collatz are sometimes known as hailstone numbers. Conway proved that the original Collatz problem has no nontrivial cycles of a positive length $<400$. Lagarias (1985) showed that there are no nontrivial cycles with length $<275 000$. Conway (1972) also proved that Collatz-type problems can be formally undecidable.
See also the work by Tao (2020) for dramatic recent progress on the problem.

The Collatz algorithm has been tested and found to always reach 1 for all numbers $\leq19\times2^{58}\approx 5.48\times10^{18}$ (Oliveira e Silva 2008), for all numbers $\leq10^{20}$, and most recently for all numbers $\leq2^{68}\approx2.95\times10^{20}$ (Tao 2020). Because of the difficulty in solving this problem, Erd\H{o}s commented that "mathematics is not yet ready for such problems" (Lagarias 1985)~\cite{N,S}.

Different from former researches on the Collatz problem, we form a Collatz's Sequence
\[
S_{c}(m,t)=\{n_{0}(m,t),n_{1}(m,t),n_{2}(m,t),\ldots,n_{m}(m,t),n_{m+1}(m,t)\}
\]
produced by a function $n_{i+1}(m,t)=(3n_{i}(m,t)+1)/2$ for $i=0,1,2,\ldots,m$ and ended by an even number $n_{m+1}(m,t)$(See Definition~\ref{def:odd_seq_1}), and prove
\[
n_{i}(m,t)=2^{m+1-i}\times3^{i}t-1\texttt{ for }i=0,1,2,\ldots,m+1,
\]
by which all odd numbers can be expressed. Then we discuss and solve the Collatz problem in two ways:

(1) In the proof of Theorem~\ref{th:Collatz_s}, by theories of the matrixes and Diophantine Equations we prove that no nontrivial Collatz cycle(See Definition~\ref{def:odd_seq_c}) exists in a series of Collatz's Sequences(See Definition~\ref{def:odd_seq_s})
\[
S_{c}(s)=\{S_{c}(m_{j},t_{j})|j=0,1,2,\ldots,s-1\}\texttt{ where }0<s<\infty,
\]
then we derive that the odd number $n_{i_{s}}(m_{s},t_{s})$ in a series of Collatz's Sequences $S_{c}(s)$ is equal to 1 for $s$ finite, and finally we prove that the odd number $n_{i_{s}}(m_{s},t_{s})$ in a series of Collatz's Sequences $S_{c}(s)$ is equal to 1 as $s$ tends to infinite based on the distribution of even numbers and the law of great numbers.

(2) In the proof of Theorem~\ref{th:Collatz_1}, we prove that first all odd numbers in the range $(N,\rho N]$(See Definition~\ref{def:odd_seq_rho}) return to 1 where these lemmas~\ref{le:odd_seqs_n_m_m-1}-\ref{le:odd_seqs_n_m_m-1_r1} form mutual references or nested recursions to themselves, which make complexities of the Collatz problem, and next by letting $N_{0}=N$ and $N_{i+1}=\rho N_{i}$ each odd number in the new range $(N_{i+1},\rho N_{i+1}]$ returns to 1 for $i=0,1,2,\ldots$. So we prove that all odd numbers return to 1.

Hence the solution of the Collatz problem is derived, which is stated as a theorem: Each Collatz's Sequence always returns to 1.

\section{Definitions of Collatz's Sequences}

\begin{definition}\label{def:odd_seq_1}
Given a nonnegative integer $m$ and an odd number $t$ not divisible by $3$, let a Collatz's Sequence be
\[
S_{c}(m,t)=\{n_{0}(m,t),n_{1}(m,t),n_{2}(m,t),\ldots,n_{m}(m,t),n_{m+1}(m,t)\}
\]
produced by a function $n_{i+1}(m,t)=(3n_{i}(m,t)+1)/2$ for $i=0,1,2,\ldots,m$ and ended by an even number $n_{m+1}(m,t)$ where $n_{0}(m,t)$ is odd and no odd number $n$ exists such that $n_{0}(m,t)=(3n+1)/2$. Then each odd number $n$ can be expressed as an odd number $n_{i}(m,t)\in S_{c}(m,t)$ where $0\leq i\leq m$, i.e., each odd number $n$ implies a Collatz's Sequence $S_{c}(m,t)$.
\end{definition}

Then given a Collatz's Sequence $S_{c}(m,t)$, when it returns to 1 each number $n_{i}(m,t)\in S_{c}(m,t)$ for $i=0,1,2,\ldots,m,m+1$ returns to 1, and when any number $n_{i}(m,t)\in S_{c}(m,t)$ for $i=0,1,2,\ldots,m,m+1$ returns to 1 it returns to 1.

\begin{definition}\label{def:odd_seq_s}
Given a nonnegative integer $m$ and an odd number $t$ not divisible by $3$, let a series of Collatz's Sequences be
\[
S_{c}(s)=\{S_{c}(m_{j},t_{j})|j=0,1,2,\ldots,s-1\}
\]
where $m_{0}=m$, $t_{0}=t$, and for $j=0,1,2,\ldots,s-1$ the odd number $n_{i_{j}}(m_{j},t_{j})$ is greater than 1, the odd number $n_{i_{j+1}}(m_{j+1},t_{j+1})$ is not divisible by $3$, and the even number
\[
n_{m_{j}+1}(m_{j},t_{j})=2^{r_{j+1}}n_{i_{j+1}}(m_{j+1},t_{j+1})\texttt{ where }r_{j+1}\geq1\texttt{ and }0\leq i_{j+1}\leq m_{j+1}.
\]
\end{definition}

\begin{definition}\label{def:odd_seq_c}
Let a Collatz cycle be a cycle starting from an odd number $n$ and returning to $n$ by the Collatz algorithm. If a Collatz cycle contains at least one odd number greater than 1 then it is nontrivial else it is trivial.
\end{definition}

\begin{definition}\label{def:odd_seq_rc}
Let define
\[
r_{c}(m,t,i,1)=\frac{n_{i}(m,t)}{n_{i_{1}}(m_{1},t_{1})}\texttt{ and }r_{c}(m,t,i,s)=\frac{n_{i}(m,t)}{n_{i_{s}}(m_{s},t_{s})}
\]
where $0\leq i\leq m$, $0\leq i_{1}\leq m_{1}$, and $0\leq i_{s}\leq m_{s}$.
\end{definition}

\begin{definition}\label{def:odd_seq_rc_j}
Let define
\[
r_{c}(m_{j},t_{j})=\frac{n_{i_{j}}(m_{j},t_{j})}{n_{i_{j+1}}(m_{j+1},t_{j+1})}\texttt{ for }j=0,1,2,\ldots,s-1.
\]
Then there is
\[
r_{c}(m,t,i,s)=\prod_{j=0}^{s-1}r_{c}(m_{j},t_{j})\texttt{ where }i_{0}=i,\texttt{ }m_{0}=m,\texttt{ and }t_{0}=t.
\]
\end{definition}

\begin{definition}\label{def:odd_seq_n}
Since the Collatz algorithm has been tested and found to always reach 1 for all numbers $\leq2^{68}$ (Tao 2020), let the Collatz algorithm have been tested and found to always reach 1 for all numbers less than or equal to a nature number $N$. For example, the first nature number $N$ is equal to $2^{68}$.
\end{definition}

\begin{definition}\label{def:odd_seq_rho}
Let define a current nature number range $(N,\rho N]$ where $\rho$ is a rational number satisfying $1<\rho=\rho(\kappa)\leq\rho(3)$ where $\rho(\kappa)=1+1/2^{\kappa}$ for $\kappa\geq3$. We let $\rho=\rho(3)$.
\end{definition}

\begin{definition}\label{def:sigma_eta}
For $\xi=1,2,\cdots,s$, let define $\eta=\eta_{s}$, $\eta_{0}=0$,
\[
\sigma(\zeta)=\sum_{j=1}^{\zeta}r_{j}\texttt{, }\sigma(\xi,\zeta)=\sum_{j=\xi}^{\zeta}r_{j}\texttt{, and }\eta_{\xi}=\sum_{j=0}^{\xi-1}(m_{j}-i_{j}+1)\texttt{ where }0\leq i_{j}\leq m_{j}.
\]
\end{definition}

\begin{definition}\label{def:matrix_a}
Let define a $\eta\times\eta$ nonsingular matrix
\[
\mathbf{A}=
\left(
  \begin{array}{ccccccccccccc}
    3 & -2 & 0 & \cdots  & 0 & 0 & 0 & \cdots  & 0 & 0 & 0 & \cdots & 0 \\
    0 & 3 & -2 & \cdots  & 0 & 0 & 0 & \cdots  & 0 & 0 & 0 & \cdots & 0 \\
    0 & 0 & 3 & \cdots  & 0 & 0 & 0 & \cdots  & 0 & 0 & 0 & \cdots & 0 \\
    \vdots & \vdots & \vdots & \cdots & \vdots & \vdots & \vdots & \cdots & \vdots & \vdots & \vdots & \cdots & \vdots \\
    0 & 0 & 0 & \cdots & 3 & -2^{1+r_{1}} & 0 & \cdots  & 0 & 0 & 0 & \cdots & 0 \\
    0 & 0 & 0 & \cdots & 0 & 3 & -2 & \cdots  & 0 & 0 & 0 & \cdots & 0 \\
    0 & 0 & 0 & \cdots & 0 & 0 & 3 & \cdots  & 0 & 0 & 0 & \cdots & 0 \\
    \vdots & \vdots & \vdots & \cdots & \vdots & \vdots & \vdots & \cdots & \vdots & \vdots & \vdots & \cdots & \vdots \\
    0 & 0 & 0 & \cdots & 0 & 0 & 0 & \cdots & 3 & -2^{1+r_{2}} & 0 & \cdots & 0 \\
    0 & 0 & 0 & \cdots & 0 & 0 & 0 & \cdots & 0 & 3 & -2 & \cdots & 0 \\
    0 & 0 & 0 & \cdots & 0 & 0 & 0 & \cdots & 0 & 0 & 3 & \cdots & 0 \\
    \vdots & \vdots & \vdots & \cdots & \vdots & \vdots & \vdots & \cdots & \vdots & \vdots & \vdots & \cdots & \vdots \\
    -2^{1+r_{s}} & 0 & 0 & \cdots & 0 & 0 & 0 & \cdots & 0 & 0 & 0 & \cdots & 3 \\
 \end{array}
\right)
\]
\[
=
\left(
  \begin{array}{cccc}
    a_{1,1} & a_{1,2} & \cdots & a_{1,\eta} \\
    a_{2,1} & a_{2,2} & \cdots & a_{2,\eta} \\
    \vdots & \vdots & \cdots & \vdots \\
    a_{\eta,1} & a_{\eta,2} & \cdots & a_{\eta,\eta} \\
 \end{array}
\right)\texttt{. Then }
\mathbf{A}^{*}=
\left(
  \begin{array}{cccc}
    A_{1,1} & A_{2,1} & \cdots & A_{\eta,1} \\
    A_{1,2} & A_{2,2} & \cdots & A_{\eta,2} \\
    \vdots & \vdots & \cdots & \vdots \\
    A_{1,\eta} & A_{2,\eta} & \cdots & A_{\eta,\eta} \\
 \end{array}
\right)
\]
is the adjoint matrix of $\mathbf{A}$ and the determinant of $\mathbf{A}$ is $Det\mathbf{A}=3^{\eta}-2^{\eta+\sigma(s)}$ where $A_{i,j}$ is the algebraic complement of the element $a_{i,j}$ in the determinant $Det\mathbf{A}$.
\end{definition}

\section{Theorems of Collatz's Sequences}

\begin{theorem}\label{th:Collatz_s}
Let $m_{0}^{*}$ be an integer. Then iterating
\begin{equation}\label{def:collatz_orig_2}
m_{k}^{*}=\{\begin{array}{cc}
         m_{k-1}^{*}/2 & \texttt{ for }m_{k-1}^{*}\texttt{ even } \\
        3m_{k-1}^{*}+1 & \texttt{ for }m_{k-1}^{*}\texttt{ odd }
      \end{array}
\end{equation}
always returns to 1 for any positive $m_{0}^{*}$, i.e., each Collatz's Sequence always returns to 1.
\end{theorem}

\begin{proof}
First, by Lemma~\ref{le:odd_seq}, given a Collatz's Sequence
\[
S_{c}(m,t)=\{n_{0}(m,t),n_{1}(m,t),n_{2}(m,t),\ldots,n_{m}(m,t),n_{m+1}(m,t)\},
\]
there are $n_{i}(m,t)=2^{m+1-i}\times3^{i}t-1$ for $i=0,1,2,\ldots,m,m+1$, and
\[
n_{m+1}(m,t)=(3n_{m}(m,t)+1)/2=3^{m+1}t-1=2^{r_{1}}n_{i_{1}}(m_{1},t_{1})
\]
where $r_{1}\geq1$, $0\leq i_{1}\leq m_{1}$, and $n_{i_{1}}(m_{1},t_{1})$ is odd and not divisible by $3$.

Second, by Lemma~\ref{le:odd_seq_s_c}, no nontrivial Collatz cycle exists in any series of Collatz's Sequences $S_{c}(s)$ where $0<s<\infty$.

Third, by Lemma~\ref{le:odd_seqs}, given a Collatz's Sequence $S_{c}(m,t)$, there exist a unique series of Collatz's Sequences $S_{c}(s)=\{S_{c}(m_{j},t_{j})|j=0,1,2,\ldots,s-1\}$ where $m_{0}=m$, $t_{0}=t$, the odd number $n_{i_{j}}(m_{j},t_{j})$ is greater than 1, the even number
\[
n_{m_{j}+1}(m_{j},t_{j})=2^{r_{j+1}}n_{i_{j+1}}(m_{j+1},t_{j+1})\texttt{ where }r_{j+1}\geq1\texttt{ and }0\leq i_{j+1}\leq m_{j+1},
\]
for $j=0,1,2,\ldots,s-1$, and the odd number $n_{i_{s}}(m_{s},t_{s})$ is equal to 1 for $s$ finite.

Finally, by Lemma~\ref{le:odd_seq_s_e}, no nontrivial Collatz cycle exists in a series of Collatz's Sequences $S_{c}(s)=\{S_{c}(m_{j},t_{j})|j=0,1,2,\ldots,s-1\}$ and the odd number $n_{i_{s}}(m_{s},t_{s})$ in the series of Collatz's Sequences $S_{c}(s)$ is equal to 1 as $s$ tends to infinite, so that each initiate odd number $n_{0}(m,t)$ returns to 1.

So, for each Collatz's Sequence $S_{c}(m,t)$, there always exists a finite positive integer $s$ such that the odd number $n_{i_{s}}(m_{s},t_{s})=n_{0}(0,1)=1$. Then starting from the odd number $n_{0}(m,t)\in S_{c}(m,t)$ the unique series of Collatz's Sequences $S_{c}(s)$ return to 1.

Hence each odd number always returns to 1, i.e., each Collatz's Sequence always returns to 1.

This completes the proof of the theorem and the Collatz conjecture is proved.
\end{proof}

\begin{theorem}\label{th:Collatz_1}
Let $m_{0}^{*}$ be an integer. Then iterating
\begin{equation}\label{def:collatz_orig_1}
m_{k}^{*}=\{\begin{array}{cc}
         m_{k-1}^{*}/2 & \texttt{ for }m_{k-1}^{*}\texttt{ even } \\
        3m_{k-1}^{*}+1 & \texttt{ for }m_{k-1}^{*}\texttt{ odd }
      \end{array}
\end{equation}
always returns to 1 for any positive $m_{0}^{*}$, i.e., each Collatz's Sequence always returns to 1.
\end{theorem}

\begin{proof}
First, given a Collatz's Sequence
\[
S_{c}(m,t)=\{n_{0}(m,t),n_{1}(m,t),n_{2}(m,t),\ldots,n_{m}(m,t),n_{m+1}(m,t)\},
\]
by Lemma~\ref{le:odd_seq}, there are
\[
n_{i}(m,t)=2^{m+1-i}\times3^{i}t-1\texttt{ for }i=0,1,2,\ldots,m,m+1,
\]
and
\[
n_{m+1}(m,t)=(3n_{m}(m,t)+1)/2=3^{m+1}t-1=2^{r_{1}}n_{i_{1}}(m_{1},t_{1})
\]
where $r_{1}\geq1$, $0\leq i_{1}\leq m_{1}$, and $n_{i_{1}}(m_{1},t_{1})$ is odd and not divisible by $3$.

By Lemma~\ref{le:odd_seq_t} the odd number $t$ can not be divisible by $3$. By Lemma~\ref{le:odd_seq_n0} the odd number $n_{0}(m,t)$ can not be expressed as $6x-1$ nor $6x+5$.

Second, by Lemma~\ref{le:odd_seq_n0_4k+1}, each odd number $n_{0}(m,t)=4k+1$ in the range $(N,\rho N]$ returns to 1, so that one half of odd numbers in the range $(N,\rho N]$ return to 1.

Third, by Lemma~\ref{le:odd_seq_ni_6x-1}, each odd number $n=6x-1$ in the range $(N,\rho N]$ returns to 1, so that another $1/6$ of odd numbers in the range $(N,\rho N]$ return to 1.

Hence $2/3$ of odd numbers in the range $(N,\rho N]$ return to 1.

Fourth, by Lemma~\ref{le:odd_seq_n0_4k+3} all odd numbers $n_{0}(m,t)=4k+3$ in the range $(N,\rho N]$ return to 1, so that another $1/3$ of odd numbers in the range $(N,\rho N]$ return to 1.

Hence all odd numbers in the range $(N,\rho N]$ return to 1.

Finally, let a new nature number $N_{1}=\rho N$. Then similarly to the above analysis, each odd number in the new range $(N_{1},\rho N_{1}]$ always returns to 1. We can go on in this way so that each odd number in the range $(N_{i+1},\rho N_{i+1}]$ always returns to 1 where the nature number $N_{i+1}=\rho N_{i}$ for $i=1,2,\ldots$. So we prove that all odd numbers return to 1.

Hence each odd number always returns to 1, i.e., each Collatz's Sequence always returns to 1.

This completes the proof of the theorem and the Collatz conjecture is proved.
\end{proof}

\section{Lemmas of Collatz's Sequences}

\subsection{Formation of Collatz's Sequences}

\begin{lemma}[Key expressions of Collatz's Sequences]\label{le:odd_seq}
Given a Collatz's Sequence
\[
S_{c}(m,t)=\{n_{0}(m,t),n_{1}(m,t),n_{2}(m,t),\ldots,n_{m}(m,t),n_{m+1}(m,t)\}
\]
where $n_{0}(m,t)$ is odd, each $n_{i}(m,t)=(3n_{i-1}(m,t)+1)/2$ for $i=1,2,\ldots,m$ is odd, and the natural number $n_{m+1}(m,t)=(3n_{m}(m,t)+1)/2$ is even, there are
\[
n_{i}(m,t)=2^{m+1-i}\times3^{i}t-1\texttt{ for }i=0,1,2,\ldots,m+1
\]
where the positive integer $t$ is odd and not divisible by $3$.
\end{lemma}

\begin{proof}
Let $n_{0}(m,t)=2x+1$ and $n_{m}(m,t)=2y+1$. Then since there are
\[
n_{m}(m,t)=(3n_{m-1}(m,t)+1)/2=(3^{2}n_{m-2}+3+2)/2^{2}
\]
\[
=(3^{3}n_{m-3}+3^{2}+3\times2+2^{2})/2^{3}=\cdots
=(3^{m}n_{0}(m,t)+\sum_{j=0}^{m-1}3^{m-1-j}2^{j})/2^{m}
\]
\[
=(3^{m}n_{0}(m,t)+3^{m}-2^{m})/2^{m}
=3^{m}(n_{0}(m,t)+1)/2^{m}-1,
\]
we have $2^{m}(y+1)=3^{m}(x+1)$. So there must be $x=2^{m}t-1\texttt{ and }y=3^{m}t-1$ or
\[
n_{0}(m,t)=2^{m+1}t-1\texttt{ and }n_{m}(m,t)=2\times3^{m}t-1,
\]
so that by $n_{i}(m,t)=(3n_{i-1}(m,t)+1)/2$ for $i=1,2,\ldots,m$ there are
\[
n_{i}(m,t)=2^{m+1-i}\times3^{i}t-1\texttt{ for }i=0,1,2,\ldots,m+1
\]
where the positive integer $t$ must be odd and not divisible by $3$ since the number $n_{m+1}(m,t)$ is even and no odd number $n$ exists such that $n_{0}(m,t)=(3n+1)/2$.

Hence, when the positive integer $m$ is finite, the number of elements of the Collatz's Sequence $S_{c}(m,t)$ is finite and equal to $m+2$.

This completes the proof of the lemma.
\end{proof}

\begin{lemma}[Property of the nonnegative integer $m$]\label{le:odd_seq_m}
Given a Collatz's Sequence
\[
S_{c}(m,t)=\{n_{0}(m,t),n_{1}(m,t),n_{2}(m,t),\ldots,n_{m}(m,t),n_{m+1}(m,t)\}
\]
where $n_{i}(m,t)=2^{m+1-i}\times3^{i}t-1$ for $i=0,1,2,\ldots,m+1$, if $n_{0}(m,t)=4k+1$ then there are $m=0$ and $n_{m+1}(m,t)=2^{r}(2k_{1}+1)$ where $r$ is a positive integer, else for $n_{0}(m,t)=4k+3$ there are $m>0$, $n_{i}(m,t)=4k_{i}+3$ for $i=1,2,\ldots,m-1$, $n_{m}(m,t)=4k_{m}+1$, and $n_{m+1}(m,t)=2^{r}(2k_{m+1}+1)$ where $r\geq1$.

On the other hand, given an odd number $n_{i}(m,t)$, nonnegative integers $m$ and $i$, and an odd number $t$ can be determined by the formula
\[
n_{i}(m,t)=2^{m+1-i}\times3^{i}t-1\texttt{ or }n_{i}(m,t)+1=2^{m+1-i}\times3^{i}t
\]
where $i=0,1,2,\ldots,m$ and no odd number $n$ exists such that $n_{0}(m,t)=(3n+1)/2$.
\end{lemma}

\begin{proof}
Given a Collatz's Sequence
\[
S_{c}(m,t)=\{n_{0}(m,t),n_{1}(m,t),n_{2}(m,t),\ldots,n_{m}(m,t),n_{m+1}(m,t)\},
\]
by Lemma~\ref{le:odd_seq}, $n_{i}(m,t)=2^{m+1-i}\times3^{i}t-1$ for $i=0,1,2,\ldots,m+1$.

Since $n_{0}(m,t)=2^{m+1}t-1=2k+1$, there is an equation $2^{m}t=k+1$. If $k$ is even or $n_{0}(m,t)=4k+1$ then there are $m=0$ and $n_{m+1}(m,t)=2^{r}(2k_{1}+1)$ where $r$ is a positive integer, else for $k$ odd or $n_{0}(m,t)=4k+3$ there is $m>0$.

When $m>0$, there are $n_{i}(m,t)=2^{m+1-i}\times3^{i}t-1=2k+1$ and $2^{m-i}\times3^{i}t=k+1$ for $i=1,2,\ldots,m$. If $i<m$ then $k$ is odd so that $n_{i}(m,t)=4k_{i}+3$, else $k$ is even so that $n_{m}(m,t)=4k_{m}+1$ and $n_{m+1}(m,t)=2^{r}(2k_{m+1}+1)$ where $r\geq1$.

On the other hand, given an odd number $n_{i}(m,t)$, nonnegative integers $m$ and $i$, and an odd number $t$ can be determined by the formula
\[
n_{i}(m,t)=2^{m+1-i}\times3^{i}t-1\texttt{ or }n_{i}(m,t)+1=2^{m+1-i}\times3^{i}t
\]
where $i=0,1,2,\ldots,m$ and no odd number $n$ exists such that $n_{0}(m,t)=(3n+1)/2$.

This completes the proof of the lemma.
\end{proof}

\begin{lemma}[Property of the odd number $t$]\label{le:odd_seq_t}
Given a Collatz's Sequence
\[
S_{c}(m,t)=\{n_{0}(m,t),n_{1}(m,t),n_{2}(m,t),\ldots,n_{m}(m,t),n_{m+1}(m,t)\}
\]
where $n_{i}(m,t)=2^{m+1-i}\times3^{i}t-1$ for $i=0,1,2,\ldots,m+1$, the odd number $t$ can not be divisible by $3$.
\end{lemma}

\begin{proof}
If the odd number $t$ can be divisible by $3$, then let an odd number $t^{*}=2x+1$. Then we have
\[
n_{0}(m,t)=2^{m+1}t-1=2^{m+1}(6x+3)-1=2^{m+1}\times3t^{*}-1.
\]

On the other hand, let an odd number $n=n_{0}(m+1,t^{*})=2^{m+2}t^{*}-1$. Then we have
\[
n_{1}(m+1,t^{*})=(3n_{0}(m+1,t^{*})+1)/2=2^{m+1}\times3t^{*}-1=n_{0}(m,t).
\]

It contradicts the fact that $n_{0}(m,t)$ is an initiate odd number in the Collatz's Sequence $S_{c}(m,t)$ and no odd number $n$ exists such that $n_{0}(m,t)=(3n+1)/2$. Hence the odd number $t$ can not be divisible by $3$.

This completes the proof of the lemma.
\end{proof}

\begin{lemma}[Property of the odd number $n_{0}(m,t)$]\label{le:odd_seq_n0}
Given a Collatz's Sequence
\[
S_{c}(m,t)=\{n_{0}(m,t),n_{1}(m,t),n_{2}(m,t),\ldots,n_{m}(m,t),n_{m+1}(m,t)\}
\]
where $n_{i}(m,t)=2^{m+1-i}\times3^{i}t-1$ for $i=0,1,2,\ldots,m+1$, the odd number $n_{0}(m,t)\in S_{c}(m,t)$ can not be expressed as $6x-1$ nor $6x+5$.
\end{lemma}

\begin{proof}
By $n_{0}(m,t)=2^{m+1}t-1$, let $2^{m+1}=6x_{0}+2$ or $2^{m+1}=6x_{0}+4$.

By Lemma~\ref{le:odd_seq_t} the odd number $t$ can not be divisible by $3$.

If $t=6x_{t}+1$, then there are
\[
n_{0}(m,t)=(6x_{0}+2)(6x_{t}+1)-1=6(6x_{0}x_{t}+x_{0}+2x_{t})+2-1=6x+1
\]
or
\[
n_{0}(m,t)=(6x_{0}+4)(6x_{t}+1)-1=6(6x_{0}x_{t}+x_{0}+4x_{t})+4-1=6x+3.
\]

Otherwise for $t=6x_{t}+5$, there are
\[
n_{0}(m,t)=(6x_{0}+2)(6x_{t}+5)-1=6(6x_{0}x_{t}+5x_{0}+2x_{t})+10-1=6x+3
\]
or
\[
n_{0}(m,t)=(6x_{0}+4)(6x_{t}+5)-1=6(6x_{0}x_{t}+5x_{0}+4x_{t})+20-1=6x+1.
\]

Hence the odd number $n_{0}(m,t)$ can not be expressed as $6x-1$ nor $6x+5$.

This completes the proof of the lemma.
\end{proof}

\begin{lemma}[Property of the positive integer $r_{1}$]\label{le:odd_seq_r1}
Given a Collatz's Sequence
\[
S_{c}(m,t)=\{n_{0}(m,t),n_{1}(m,t),n_{2}(m,t),\ldots,n_{m}(m,t),n_{m+1}(m,t)\}
\]
where $n_{i}(m,t)=2^{m+1-i}\times3^{i}t-1$ for $i=0,1,2,\ldots,m+1$, let the even number
\[
n_{m+1}(m,t)=2^{r_{1}}n_{i_{1}}(m_{1},t_{1})\texttt{ where }r_{1}>0\texttt{, }0\leq i_{1}\leq m_{1},
\]
and $n_{i_{1}}(m_{1},t_{1})$ is odd and not dividable by 3.

Then $1/2^{r}$ of positive integers $r_{1}$ satisfy $r_{1}=r$ for $r=1,2,\ldots$.
\end{lemma}

\begin{proof}
First, there is $r_{1}=1$ when $m$ is even and $t=4k+1$ where $k\geq0$ since
\[
n_{m+1}(m,t)=3^{m+1}t-1=3^{m+1}4k+3^{m+1}-1
\]
\[
=2(3^{m+1}2k+\sum_{j=0}^{m}3^{j})\texttt{ where }\sum_{j=0}^{m}3^{j}\texttt{ is odd},
\]
or when $m$ is odd and $t=4k+3$ where $k>0$ since
\[
n_{m+1}(m,t)=3^{m+1}t-1=3^{m+1}4k+3^{m+2}-1
\]
\[
=2(3^{m+1}2k+\sum_{j=0}^{m+1}3^{j})\texttt{ where }\sum_{j=0}^{m+1}3^{j}\texttt{ is odd}.
\]

Second, there is $r_{1}=2$ when $m$ is even and $t=8k-1$ where $k>0$ since
\[
n_{m+1}(m,t)=3^{m+1}t-1=3^{m+1}8k-(3^{m+1}+1)
\]
\[
=2^{2}(3^{m+1}2k-\sum_{j=0}^{m}(-1)^{j}3^{j})\texttt{ where }\sum_{j=0}^{m}(-1)^{j}3^{j}\texttt{ is odd},
\]
or when $m$ is odd and $t=8k-3$ where $k>0$ since
\[
n_{m+1}(m,t)=3^{m+1}t-1=3^{m+1}8k-(3^{m+2}+1)
\]
\[
=2^{2}(3^{m+1}2k-\sum_{j=0}^{m+1}(-1)^{j}3^{j})\texttt{ where }\sum_{j=0}^{m+1}(-1)^{j}3^{j}\texttt{ is odd}.
\]

Third, there is $r_{1}=3$ when $m=4h+1-j$ and $t=16k+3^{j}$ for $j=0,1,2,3$ where $h>0$ and $k>0$ since
\[
n_{m+1}(m,t)=3^{m+1}t-1=3^{m+1}16k+3^{4h+2}-1
\]
\[
=2^{3}(3^{m+1}2k-\sum_{j=0}^{2h}3^{2j})\texttt{ where }\sum_{j=0}^{2h}3^{2j}\texttt{ is odd}.
\]

Forth, there is $r_{1}=4$ when $m=8h+3-j$ and $t=32k+3^{j}$ for $j=0,1,2,\ldots,7$ where $h>0$ and $k>0$ since
\[
n_{m+1}(m,t)=3^{m+1}t-1=3^{m+1}32k+3^{8h+4}-1
\]
\[
=2^{4}(3^{m+1}2k-5\sum_{j=0}^{2h}3^{4j})\texttt{ where }\sum_{j=0}^{2h}3^{4j}\texttt{ is odd}.
\]

Fifth, there is $r_{1}=5$ when $m=16h+7-j$ and $t=64k+3^{j}$ for $j=0,1,2,\ldots,15$ where $h>0$ and $k>0$ since
\[
n_{m+1}(m,t)=3^{m+1}t-1=3^{m+1}64k+3^{16h+8}-1
\]
\[
=2^{5}(3^{m+1}2k-205\sum_{j=0}^{2h}3^{8j})\texttt{ where }\sum_{j=0}^{2h}3^{8j}\texttt{ is odd}.
\]

Similarly, there are $r_{1}=r$ for $r=6,7,\ldots$ when $m=2^{r-1}h+2^{r-2}-1-j$ and $t=2^{r+1}k+3^{j}$ for $j=0,1,2,\ldots,2^{r-1}-1$ where $h>0$ and $k>0$ since
\[
n_{m+1}(m,t)=3^{m+1}t-1=3^{m+1}2^{r+1}k+3^{2^{r-1}h+2^{r-2}}-1
\]
\[
=2^{r}(3^{m+1}2k-\frac{3^{2^{r-2}}-1}{2^{r}}\sum_{j=0}^{2h}3^{2^{r-2}j})\texttt{ where }\sum_{j=0}^{2h}3^{2^{r-2}j}\texttt{ is odd.}
\]

Hence $1/2^{r}$ of positive integers $r_{1}$ satisfy $r_{1}=r$ for $r=1,2,\ldots$.

This completes the proof of the lemma.
\end{proof}

\subsection{Properties of Collatz's Sequences}

\begin{lemma}[On ratio $r_{c}(m,t,i,1)$]\label{le:odd_seq_1}
Given a Collatz's Sequence
\[
S_{c}(m,t)=\{n_{0}(m,t),n_{1}(m,t),n_{2}(m,t),\ldots,n_{m}(m,t),n_{m+1}(m,t)\}
\]
where $n_{i}(m,t)=2^{m+1-i}\times3^{i}t-1$ for $i=0,1,2,\ldots,m+1$, let
\[
n_{m+1}(m,t)=2^{r_{1}}n_{i_{1}}(m_{1},t_{1})\texttt{ where }r_{1}\geq1\texttt{, }0\leq i_{1}\leq m_{1},
\]
and $n_{i_{1}}(m_{1},t_{1})$ is odd and not divisible by $3$. Then for $i=0,1,2,\ldots,m$ we have
\begin{equation}\label{eq:odd_seq_1}
r_{c}(m,t,i,1)=\frac{n_{i}(m,t)}{n_{i_{1}}(m_{1},t_{1})}=(2/3)^{m+1-i}2^{r_{1}}\prod_{k=i}^{m}(1-\frac{1}{2n_{k+1}(m,t)})
\end{equation}
where $r_{1}\geq1$, $0\leq i_{1}\leq m_{1}$, and $n_{i_{1}}(m_{1},t_{1})$ is odd and not divisible by $3$.
\end{lemma}

\begin{proof}
Given a Collatz's Sequence $S_{c}(m,t)$, the even number
\[
n_{m+1}(m,t)=3^{m+1}t-1=2^{r_{1}}n_{i_{1}}(m_{1},t_{1})\texttt{ where }r_{1}\geq1\texttt{ and }0\leq i_{1}\leq m_{1}.
\]

Let an odd number $n_{i}(m,t)$ greater than 1 be in the Collatz's Sequence $S_{c}(m,t)$ where $0\leq i\leq m$. Then by the formula $n_{i+1}(m,t)=(3n_{i}(m,t)+1)/2$ for $i=0,1,2,\ldots,m$ we have
\[
n_{i}(m,t)=(2n_{i+1}(m,t)-1)/3
\]
\[
=(2/3)(1-\frac{1}{2n_{i+1}(m,t)})n_{i+1}(m,t)
\]
\[
=(2/3)^{2}\prod_{k=i}^{i+1}(1-\frac{1}{2n_{k+1}(m,t)})n_{i+2}(m,t)=\cdots
\]
\[
=(2/3)^{m+1-i}\prod_{k=i}^{m}(1-\frac{1}{2n_{k+1}(m,t)})n_{m+1}(m,t)
\]
\[
=(2/3)^{m+1-i}2^{r_{1}}\prod_{k=i}^{m}(1-\frac{1}{2n_{k+1}(m,t)})n_{i_{1}}(m_{1},t_{1}).
\]

Hence Formula~\ref{eq:odd_seq_1} holds.

This completes the proof of the lemma.
\end{proof}

\begin{lemma}\label{le:odd_seq_1_c}
No nontrivial Collatz cycle exists in any Collatz's Sequence $S_{c}(m,t)$.
\end{lemma}

\begin{proof}
Given a Collatz's Sequence
\[
S_{c}(m,t)=\{n_{0}(m,t),n_{1}(m,t),n_{2}(m,t),\ldots,n_{m}(m,t),n_{m+1}(m,t)\}
\]
where $n_{i}(m,t)=2^{m+1-i}\times3^{i}t-1$ for $i=0,1,2,\ldots,m+1$, the even number
\[
n_{m+1}(m,t)=3^{m+1}t-1=2^{r_{1}}n_{i_{1}}(m_{1},t_{1})\texttt{ where }r_{1}\geq1\texttt{ and }0\leq i_{1}\leq m_{1}.
\]

Let an odd number $n_{i}(m,t)$ greater than 1 be in the Collatz's Sequence $S_{c}(m,t)$ where $0\leq i\leq m$ and consider the ratio of the odd number $n_{i}(m,t)$ to the odd number $n_{i_{1}}(m_{1},t_{1})$. Then we have
\[
r_{c}(m,t,i,1)=\frac{n_{i}(m,t)}{n_{i_{1}}(m_{1},t_{1})}=\frac{2^{m+1-i}\times3^{i}t-1}{(3^{m+1}t-1)/2^{r_{1}}}.
\]

Thus if $r_{c}(m,t,i,1)=1$ then there must be
\[
(2^{m+1-i+r_{1}}-3^{m+1-i})3^{i}t=2^{r_{1}}-1\texttt{ and }
\frac{3^{m+1}t-1}{2^{m+1-i}\times3^{i}t-1}=2^{r_{1}}.
\]

Since a group of numbers $i=0$, $m=0$, $t=1$, and $r_{1}=1$ satisfy the above two equalities, then when $n_{i}(m,t)=n_{0}(0,1)=1$ and $r_{1}=1$, there is $r_{c}(m,t,i,1)=1$.

First, if $(2^{m+1-i+r_{1}}-3^{m+1-i})3^{i}t<2^{r_{1}}-1$ then there is $r_{c}(m,t,i,1)<1$.

Second, if $i=m=0$ and $t>0$ or $m-i=0$ and $m>0$ then there are
\[
r_{c}(m,t,i,1)=\frac{2^{m+1-i}\times3^{i}t-1}{(3^{m+1}t-1)/2^{r_{1}}}\geq\frac{4\times3^{m}t-2}{3\times3^{m}t-1}>1
\]
else if $r_{1}\geq m+1-i>1$ then there is $r_{c}(m,t,i,1)>1$ since for $r_{1}=m+1-i$
\[
(2^{m+1-i+r_{1}}-3^{m+1-i})3^{i}t\geq(4^{m+1-i}-3^{m+1-i})3^{i}
\]
\[
=3^{i}\sum_{j=0}^{m-i}4^{j}3^{m-i-j}>(m+1-i)3^{m}>2^{r_{1}}-1.
\]

Finally, for $r_{1}=1$, there are
\[
(2^{m+1-i+r_{1}}-3^{m+1-i})3^{i}t=(2\times2^{m+1-i}-3^{m+1-i})3^{i}t\texttt{ and }2^{r_{1}}-1=1.
\]
If $i>0$ then there are $3^{i}t>1$ and $(2\times2^{m+1-i}-3^{m+1-i})3^{i}t\neq1$, so that $r_{c}(m,t,i,1)\neq1$, else for $i=0$ only if $m=0$ and $t=1$ then there is $(2\times2^{m+1}-3^{m+1})t=1$, i.e., only if $n_{i}(m,t)=n_{0}(0,1)=1$ and $r_{1}=1$ then there is $r_{c}(m,t,i,1)=1$.

Hence no nontrivial Collatz cycle exists in $S_{c}(m,t)$.

This completes the proof of the lemma.
\end{proof}

\begin{lemma}[On ratio $r_{c}(m,t,i,s)$]\label{le:odd_seq_s}
Given a series of Collatz's Sequences
\[
S_{c}(s)=\{S_{c}(m_{j},t_{j})|j=0,1,2,\ldots,s-1\}
\]
where $m_{0}=m$, $t_{0}=t$, and for $j=0,1,2,\ldots,s-1$ the even number
\[
n_{m_{j}+1}(m_{j},t_{j})=2^{r_{j+1}}n_{i_{j+1}}(m_{j+1},t_{j+1})\texttt{ where }r_{j+1}\geq1\texttt{ and }0\leq i_{j+1}\leq m_{j+1},
\]
then with $0\leq i\leq m$ and $i_{0}=i$ we have
\begin{equation}\label{eq:odd_seq_s}
r_{c}(m,t,i,s)=\frac{n_{i}(m,t)}{n_{i_{s}}(m_{s},t_{s})}=\prod_{j=0}^{s-1}(2/3)^{m_{j}+1-i_{j}}2^{r_{j+1}}\prod_{k=i_{j}}^{m_{j}}(1-\frac{1}{2n_{k+1}(m_{j},t_{j})})
\end{equation}
where the number $n_{i_{s}}(m_{s},t_{s})$ is odd and $0\leq i_{s}\leq m_{s}$.
\end{lemma}

\begin{proof}
Let an odd number $n_{i}(m,t)$ greater than 1 be in the first Collatz's Sequence $S_{c}(m,t)$ where $0\leq i\leq m$. Then by Lemma~\ref{le:odd_seq_1}, we have
\[
r_{c}(m,t,i,1)=\frac{n_{i}(m,t)}{n_{i_{1}}(m_{1},t_{1})}=(2/3)^{m+1-i}2^{r_{1}}\prod_{k=i}^{m}(1-\frac{1}{2n_{k+1}(m,t)})\texttt{ where }r_{1}\geq1.
\]

Similar to the analysis in the proof of Lemma~\ref{le:odd_seq_1} for $i=0,1,2,\ldots,m$ we have
\[
n_{i}(m,t)=(2/3)^{m+1-i}2^{r_{1}}\prod_{k=i}^{m}(1-\frac{1}{2n_{k+1}(m,t)})n_{i_{1}}(m_{1},t_{1})
\]
\[
=\prod_{j=0}^{1}(2/3)^{m_{j}+1-i_{j}}2^{r_{j+1}}\prod_{k=i_{j}}^{m_{j}}(1-\frac{1}{2n_{k+1}(m_{j},t_{j})})n_{i_{2}}(m_{2},t_{2})
\]
\[
=\prod_{j=0}^{2}(2/3)^{m_{j}+1-i_{j}}2^{r_{j+1}}\prod_{k=i_{j}}^{m_{j}}(1-\frac{1}{2n_{k+1}(m_{j},t_{j})})n_{i_{3}}(m_{3},t_{3})
\]
\[
=\cdots=\prod_{j=0}^{s-1}(2/3)^{m_{j}+1-i_{j}}2^{r_{j+1}}\prod_{k=i_{j}}^{m_{j}}(1-\frac{1}{2n_{k+1}(m_{j},t_{j})})n_{i_{s}}(m_{s},t_{s}).
\]

Hence Formula~\ref{eq:odd_seq_s} holds.

This completes the proof of the lemma.
\end{proof}

\begin{lemma}\label{le:odd_seq_s_c}
No nontrivial Collatz cycle exists in any series of Collatz's Sequences $S_{c}(s)$ where $0<s<\infty$.
\end{lemma}

\begin{proof}
If $s=1$ then by Lemma~\ref{le:odd_seq_1_c} no nontrivial Collatz cycle exists in any Collatz's Sequence $S_{c}(m,t)$.

Otherwise for $1<s<\infty$ given a series of Collatz's Sequences
\[
S_{c}(s)=\{S_{c}(m_{j},t_{j})|j=0,1,2,\ldots,s-1\}
\]
where $m_{0}=m$, $t_{0}=t$, and for $j=0,1,2,\ldots,s-1$ the even number
\[
n_{m_{j}+1}(m_{j},t_{j})=2^{r_{j+1}}n_{i_{j+1}}(m_{j+1},t_{j+1})\texttt{ where }r_{j+1}\geq1\texttt{ and }0\leq i_{j+1}\leq m_{j+1}.
\]

Let an odd number $n_{i}(m,t)$ greater than 1 be in the first Collatz's Sequence $S_{c}(m,t)$ where $0\leq i\leq m$, and by Definition~\ref{def:sigma_eta} there are $\eta=\eta_{s}$, $\eta_{0}=0$,
\[
\sigma(\zeta)=\sum_{j=1}^{\zeta}r_{j}\texttt{, }\sigma(\xi,\zeta)=\sum_{j=\xi}^{\zeta}r_{j}\texttt{, and } \eta_{\xi}=\sum_{j=0}^{\xi-1}(m_{j}-i_{j}+1)\texttt{ where }0\leq i_{j}\leq m_{j}
\]
for $\xi=1,2,\cdots,s$. Then by Lemma~\ref{le:odd_seq_s}, with $0\leq i\leq m$ and $i_{0}=i$ we have
\[
r_{c}(m,t,i,s)=\frac{n_{i}(m,t)}{n_{i_{s}}(m_{s},t_{s})}<\prod_{j=0}^{s-1}(2/3)^{m_{j}+1-i_{j}}2^{r_{j+1}}
\]
where the number $n_{i_{s}}(m_{s},t_{s})$ is odd and $0\leq i_{s}\leq m_{s}$.

Thus if $r_{c}(m,t,i,s)\neq1$ then no nontrivial Collatz cycle exists in $S_{c}(s)$
else if $r_{c}(m,t,i,s)=1$ then there must be $n_{i_{s}}(m_{s},t_{s})=n_{i}(m,t)$,
\[
\prod_{j=0}^{s-1}3^{m_{j}+1-i_{j}}<\prod_{j=0}^{s-1}2^{m_{j}+1-i_{j}}2^{r_{j+1}}\texttt{ or }
3^{\eta}<2^{\eta+\sigma(s)}\texttt{ where }\eta=\eta_{s},
\]
and there must exist at least one Collatz's Sequence $S_{c}(m_{h},t_{h})\in S_{c}(s)$ where $0\leq h<s$, in which the element $n_{i_{h}}(m_{h},t_{h})$ satisfies $m_{h}-i_{h}+1\geq2$, i.e., there must be $\eta>s$, otherwise there shall be $n_{i_{s}}(m_{s},t_{s})<n_{i}(m,t)$ such that $r_{c}(m,t,i,s)>1$.

By the assumption $r_{c}(m,t,i,s)=1$ for $1<s<\infty$ and by Definition~\ref{def:odd_seq_1} for $j=0,1,2,\ldots,s-1$ there are $\eta$ linear equations
\[
2n_{k+1}(m_{j},t_{j})-3n_{k}(m_{j},t_{j})=1\texttt{ for }k=i_{j},i_{j}+1,i_{j}+2,\ldots,m_{j}
\]
where $n_{m_{j}+1}(m_{j},t_{j})=2^{r_{j+1}}n_{i_{j+1}}(m_{j+1},t_{j+1})$ and $n_{i_{s}}(m_{s},t_{s})=n_{i}(m,t)$.

Let rewrite these linear equations in the form of a matrix equation. Then there is $\mathbf{A}\mathbf{X}=\mathbf{B}$ where $\mathbf{A}$ is defined by Definition~\ref{def:matrix_a} with two dimensions $\eta\times\eta$,
\[
\mathbf{B}=-
\left(
  \begin{array}{cccc}
   1 & 1 & \cdots & 1 \\
 \end{array}
\right)^{T}\texttt{ and }\mathbf{X}=
\left(
  \begin{array}{cccc}
   x_{1} & x_{2} & \cdots & x_{\eta} \\
 \end{array}
\right)^{T}
\]
are vectors with a dimension $\eta$ where for $j=0,1,2,\ldots,s-1$
\[
x_{\eta_{j}+1}=n_{i_{j}}(m_{j},t_{j})\texttt{, }x_{\eta_{j}+2}=n_{i_{j}+1}(m_{j},t_{j})\texttt{, }\ldots\texttt{ }x_{\eta_{j+1}}=n_{m_{j}}(m_{j},t_{j}).
\]

Thus for $1<s<\infty$ the determinant $Det\mathbf{A}$ of the matrix $\mathbf{A}$ satisfies
\begin{equation}\label{fm:odd_seq_Det_A}
Det\mathbf{A}=3^{\eta}-2^{\eta+\sigma(s)}<0\texttt{, }-Det\mathbf{A}=2^{\eta+\sigma(s)}-3^{\eta}\geq5\texttt{, and }\eta>s
\end{equation}
where it is obvious that $-Det\mathbf{A}$ can not be dividable by $3$ so that $-Det\mathbf{A}\neq3$, and when $\eta>1$ or $\eta+\sigma(s)>2$ if $-Det\mathbf{A}=1$, i.e., if $4(2^{\eta+\sigma(s)-2}-1)=3(3^{\eta-1}-1)$, then first for $\eta=2k+2$ where $k\geq0$, there should be $2(2^{\eta+\sigma(s)-2}-1)=3\sum_{i=0}^{2k}3^{i}$, but this equality can not hold since its left-hand value is even and its right-hand value is odd, so that $-Det\mathbf{A}\neq1$, and next for $\eta=2k+3$ where $k\geq0$, there should be $2^{\eta+\sigma(s)-2}-1=6\sum_{i=0}^{k}3^{2i}$, but this equality can not hold since its left-hand value is odd and its right-hand value is even, so that $-Det\mathbf{A}\neq1$.

Then the matrix $\mathbf{A}$ is nonsingular and there is $\mathbf{A}^{-1}=\mathbf{A}^{*}/Det\mathbf{A}$ where $\mathbf{A}^{-1}$ and $\mathbf{A}^{*}$ are the inverse matrix and the adjoint matrix of the matrix $\mathbf{A}$, respectively.

So we have the solution $\mathbf{X}$ of the matrix equation $\mathbf{A}\mathbf{X}=\mathbf{B}$ by
\[
\mathbf{X}=\mathbf{A}^{-1}\mathbf{B}=\mathbf{A}^{*}\mathbf{B}/Det\mathbf{A}
\]
and for $\kappa=1,2,\ldots,\eta$ there must be different positive odd integers
\begin{equation}\label{fm:kappa_rows_x}
x_{\kappa}=-\sum_{j=1}^{\eta}A_{j,\kappa}/Det\mathbf{A}=-A_{\kappa}/Det\mathbf{A}\texttt{ where }A_{\kappa}=\sum_{j=1}^{\eta}A_{j,\kappa}.
\end{equation}

Thus for $\kappa=1,2,\ldots,\eta$ each number $-A_{\kappa}/Det\mathbf{A}$ must be an exact division or each $A_{\kappa}$ must be exactly dividable by $-Det\mathbf{A}$.

By Lemma~\ref{le:matrix_adj_a} the adjoint matrix $\mathbf{A}^{*}$ of the matrix $\mathbf{A}$ can be determined, so that numbers $A_{\kappa}$ for $\kappa=1,2,\ldots,\eta$ can also be determined. Particularly we have
\[
A_{1}=\sum_{\zeta=1}^{s}2^{\eta_{\zeta-1}+\sigma(\zeta-1)}3^{\eta-\eta_{\zeta}}(3^{\eta_{\zeta}-\eta_{\zeta-1}}-2^{\eta_{\zeta}-\eta_{\zeta-1}})
\]
\[
=\sum_{\zeta=1}^{s}2^{\eta_{\zeta}+\sigma(\zeta-1)}3^{\eta-\eta_{\zeta}}(2^{\sigma(\zeta,\zeta)}-1)+Det\mathbf{A}.
\]

By Lemma~\ref{le:diophantine_eq_sol} or by Lemma~\ref{le:matrix_eq_sol} $A_{1}$ can not be exactly dividable by $-Det\mathbf{A}$ when $0<s<\infty$. Thus there exits no positive odd integer solution of the matrix equation $\mathbf{A}\mathbf{X}=\mathbf{B}$. It implies that the assumption $r_{c}(m,t,i,s)=1$ does not hold. Then there must be $r_{c}(m,t,i,s)\neq1$ when $0<s<\infty$.

Hence no nontrivial Collatz cycle exists in $S_{c}(s)$ where $0<s<\infty$.

This completes the proof of the lemma.
\end{proof}

\begin{lemma}\label{le:matrix_adj_a}
The adjoint matrix $\mathbf{A}^{*}$ of a $\eta\times\eta$ nonsingular matrix $\mathbf{A}$ defined by Definition~\ref{def:matrix_a} can be determined where $Det\mathbf{A}=3^{\eta}-2^{\eta+\sigma(s)}$, $\eta=\eta_{s}$, $\eta_{0}=0$,
\[
\sigma(\zeta)=\sum_{j=1}^{\zeta}r_{j}\texttt{, }\sigma(\xi,\zeta)=\sum_{j=\xi}^{\zeta}r_{j}\texttt{, and }\eta_{\xi}=\sum_{j=0}^{\xi-1}(m_{j}-i_{j}+1)\texttt{ where }0\leq i_{j}\leq m_{j}
\]
for $\xi=1,2,\cdots,s$.

Let $A_{\kappa}=\sum_{j=1}^{\eta}A_{j,\kappa}$ for $\kappa=1,2,\ldots,\eta$. Then there are
\begin{equation}\label{fm:kappa_rows_l}
A_{\kappa}=2^{\eta-\kappa+1+\sigma(s)}3^{\kappa-1}[1-(2/3)^{\kappa-1}]+3^{\eta}[1-(2/3)^{\eta_{1}-\kappa+1}]
\end{equation}
\[
+\sum_{\zeta=1}^{s-1}2^{\eta_{\zeta}+1-\kappa+\sigma(\zeta)}3^{\eta-\eta_{\zeta}+\kappa-1}[1-(2/3)^{\eta_{\zeta+1}-\eta_{\zeta}}]
\texttt{ for }\kappa=1,2,\ldots,\eta_{1}
\]
and for $\xi=1,2,\ldots,s-1$ and $\kappa=\eta_{\xi}+1,\eta_{\xi}+2,\ldots,\min(\eta_{\xi+1},\eta)$
\begin{equation}\label{fm:kappa_rows_2}
A_{\kappa}=\sum_{\zeta=0}^{\xi}2^{\eta-\kappa+\eta_{\zeta}+1+\sigma(\xi+1,s)+\sigma(\zeta)}3^{\kappa-\eta_{\zeta}-1}[1-(2/3)^{\min(\eta_{\zeta+1},\kappa-1)-\eta_{\zeta}}]
\end{equation}
\[
+3^{\eta}[1-(2/3)^{\eta_{\xi+1}-\kappa+1}]
+\sum_{\zeta=\xi+1}^{s-1}2^{\eta_{\zeta}+1-\kappa+\sigma(\xi+1,\zeta)}3^{\eta-\eta_{\zeta}+\kappa-1}[1-(2/3)^{\eta_{\zeta+1}-\eta_{\zeta}}].
\]

Particularly we have
\begin{equation}\label{fm:kappa_row_l}
A_{1}=\sum_{\zeta=1}^{s}2^{\eta_{\zeta-1}+\sigma(\zeta-1)}3^{\eta-\eta_{\zeta}}(3^{\eta_{\zeta}-\eta_{\zeta-1}}-2^{\eta_{\zeta}-\eta_{\zeta-1}})
\end{equation}
\[
=\sum_{\zeta=1}^{s}2^{\eta_{\zeta}+\sigma(\zeta-1)}3^{\eta-\eta_{\zeta}}(2^{\sigma(\zeta,\zeta)}-1)+Det\mathbf{A}.
\]
\end{lemma}

\begin{proof}
Obviously there are $A_{k,k}=3^{\eta-1}$ for $k=1,2,\ldots,\eta$.

On the one hand, we use the formula $\mathbf{A}\mathbf{A}^{*}=\mathbf{I}Det\mathbf{A}$ where $\mathbf{I}$ is a unit matrix to determine the first collum of the matrix $\mathbf{A}^{*}$. By the formula $3A_{1,1}-2A_{1,2}=Det\mathbf{A}$ there are
\[
A_{1,2}=(3A_{1,1}-Det\mathbf{A})/2=2^{\eta-1+\sigma(s)},
\]
and by the formula $3A_{1,\kappa-1}-2A_{1,\kappa}=0$ for $\kappa=3,4,\ldots,\eta_{1}$ there are
\[
A_{1,\kappa}=3A_{1,\kappa-1}/2=2^{\eta-\kappa+1+\sigma(s)}3^{\kappa-2}\texttt{ for }\kappa=3,4,\ldots,\eta_{1}.
\]
So we obtain
\begin{equation}\label{fm:first_collum_1}
A_{1,\kappa}=2^{\eta-\kappa+1+\sigma(s)}3^{\kappa-2}\texttt{ for }\kappa=2,3,\ldots,\eta_{1}.
\end{equation}

For $\xi=1,2,\ldots,s-1$ and $\kappa=\eta_{\xi}+1,\eta_{\xi}+2,\ldots,\eta_{\xi+1}$, if $\kappa=\eta_{\xi}+1$ then by the formula $3A_{1,\kappa-1}-2^{1+r_{\xi}}A_{1,\kappa}=0$, there are
\[
A_{1,\kappa}=3A_{1,\kappa-1}/2^{1+r_{\xi}}=2^{\eta-\kappa+1+\sigma(\xi+1,s)}3^{\kappa-2},
\]
else by the formula $3A_{1,\kappa-1}-2A_{1,\kappa}=0$ for $\kappa=\eta_{\xi}+2,\eta_{\xi}+3,\ldots,\eta_{\xi+1}$ there are
\[
A_{1,\kappa}=3A_{1,\kappa-1}/2=2^{\eta-\kappa+1+\sigma(\xi+1,s)}3^{\kappa-2}\texttt{ for }\kappa=\eta_{\xi}+2,\eta_{\xi}+3,\ldots,\eta_{\xi+1}.
\]
So for $\xi=1,2,\ldots,s-1$ we obtain
\begin{equation}\label{fm:first_collum_2}
A_{1,\kappa}=2^{\eta-\kappa+1+\sigma(\xi+1,s)}3^{\kappa-2}\texttt{ for }\kappa=\eta_{\xi}+1,\eta_{\xi}+2,\ldots,\eta_{\xi+1}.
\end{equation}

Thus by formulas~\ref{fm:first_collum_1}-\ref{fm:first_collum_2} the first collum of the matrix $\mathbf{A}^{*}$ is determined.

On the other hand, we use the formula $\mathbf{A}^{*}\mathbf{A}=\mathbf{I}Det\mathbf{A}$ where $\mathbf{I}$ is a unit matrix to determine the last collum and all rows of the matrix $\mathbf{A}^{*}$.

First, we determine the last collum of the matrix $\mathbf{A}^{*}$. By the formula
\[
3A_{1,1}-2^{1+r_{s}}A_{\eta,1}=Det\mathbf{A}
\]
there are
\[
A_{\eta,1}=(3A_{1,1}-Det\mathbf{A})/2^{1+r_{s}}=2^{\eta-1+\sigma(s-1)},
\]
and by the formula $3A_{1,\kappa}-2^{1+r_{s}}A_{\eta,\kappa}=0$ for $\kappa=2,3,\ldots,\eta_{1}$ there are
\[
A_{\eta,\kappa}=3A_{1,\kappa}/2^{1+r_{s}}=2^{\eta-\kappa+\sigma(s-1)}3^{\kappa-1}\texttt{ for }\kappa=2,3,\ldots,\eta_{1}.
\]
So we obtain
\begin{equation}\label{fm:last_collum_1}
A_{\eta,\kappa}=2^{\eta-\kappa+\sigma(s-1)}3^{\kappa-1}\texttt{ for }\kappa=1,2,\ldots,\eta_{1}.
\end{equation}

For $\xi=1,2,\ldots,s-1$, by the formula $3A_{1,\kappa}-2^{1+r_{s}}A_{\eta,\kappa}=0$ there are
\[
A_{\eta,\kappa}=3A_{1,\kappa}/2^{1+r_{s}}=2^{\eta-\kappa+\sigma(\xi+1,s-1)}3^{\kappa-1}
\]
for $\kappa=\eta_{\xi}+1,\eta_{\xi}+2,\ldots,\min(\eta_{\xi+1},\eta-1)$. So for $\xi=1,2,\ldots,s-1$ we obtain
\begin{equation}\label{fm:last_collum_2}
A_{\eta,\kappa}=2^{\eta-\kappa+\sigma(\xi+1,s-1)}3^{\kappa-1}\texttt{ for }\kappa=\eta_{\xi}+1,\eta_{\xi}+2,\ldots,\min(\eta_{\xi+1},\eta-1).
\end{equation}

Thus by formulas~\ref{fm:last_collum_1}-\ref{fm:last_collum_2} the last collum of the matrix $\mathbf{A}^{*}$ is determined.

Second, we determine the $\kappa^{th}$ row of the matrix $\mathbf{A}^{*}$ where $\kappa=1,2,\ldots,\eta-1$ and positive integers $A_{k,\kappa}$ for $k=\kappa+1,\kappa+2,\ldots,\eta$.

For $\xi=0,1,2,\ldots,s-1$ and $\kappa=\eta_{\xi}+1,\eta_{\xi}+2,\ldots,\min(\eta_{\xi+1},\eta-1)$, by the formula $3A_{k,\kappa}-2A_{k-1,\kappa}=0$ for $k=\kappa+1,\kappa+2,\ldots,\eta_{\xi+1}$ there are
\[
A_{k,\kappa}=2A_{k-1,\kappa}/3=2^{k-\kappa}3^{\eta-k+\kappa-1}\texttt{ for }k=\kappa+1,\kappa+2,\ldots,\eta_{\xi+1},
\]
so that we obtain
\begin{equation}\label{fm:rows_r_1}
A_{k,\kappa}=2^{k-\kappa}3^{\eta-k+\kappa-1}\texttt{ for }k=\kappa+1,\kappa+2,\ldots,\eta_{\xi+1},
\end{equation}
and for $\zeta=\xi+1,\xi+2,\ldots,s-1$ and $k=\eta_{\zeta}+1,\eta_{\zeta}+2,\ldots,\eta_{\zeta+1}$, if $k=\eta_{\zeta}+1$ then by the formula $3A_{k,\kappa}-2^{1+r_{\zeta+1}}A_{k-1,\kappa}=0$ there are
\[
A_{k,\kappa}=2^{1+r_{\zeta+1}}A_{k-1,\kappa}/3=2^{k-\kappa+\sigma(\xi+1,\zeta)}3^{\eta-k+\kappa-1},
\]
else by the formula $3A_{k,\kappa}-2A_{k-1,\kappa}=0$ there are
\[
A_{k,\kappa}=2A_{k-1,\kappa}/3
=2^{k-\kappa+\sigma(\xi+1,\zeta)}3^{\eta-k+\kappa-1}\texttt{ for }k=\eta_{\zeta}+2,\eta_{\zeta}+3,\ldots,\eta_{\zeta+1},
\]
so that for $\zeta=\xi+1,\xi+2,\ldots,s-1$ we obtain
\begin{equation}\label{fm:rows_r_2}
A_{k,\kappa}=2^{k-\kappa+\sigma(\xi+1,\zeta)}3^{\eta-k+\kappa-1}\texttt{ for }k=\eta_{\zeta}+1,\eta_{\zeta}+2,\ldots,\eta_{\zeta+1}.
\end{equation}

Thus by formulas~\ref{fm:rows_r_1}-\ref{fm:rows_r_2} integers $A_{k,\kappa}$ are determined for $\kappa=1,2,\ldots,\eta-1$ and $k=\kappa+1,\kappa+2,\ldots,\eta$.

Finally, we determine the $\kappa^{th}$ row of the matrix $\mathbf{A}^{*}$ where $\kappa=3,4,\ldots,\eta$ and positive integers $A_{k,\kappa}$ for $k=2,3,\ldots,\kappa-1$. By the formula $3A_{k,\kappa}-2A_{k-1,\kappa}=0$ for $\kappa=3,4,\ldots,\eta_{1}$ there are
\[
A_{k,\kappa}=2A_{k-1,\kappa}/3=2^{\eta-\kappa+k+\sigma(s)}3^{\kappa-k-1}\texttt{ for }k=2,3,\ldots,\kappa-1.
\]
So we obtain
\begin{equation}\label{fm:rows_l_1}
A_{k,\kappa}=2^{\eta-\kappa+k+\sigma(s)}3^{\kappa-k-1}\texttt{ for }k=1,2,\ldots,\kappa-1.
\end{equation}

For $\xi=1,2,\ldots,s-1$ and $\kappa=\eta_{\xi}+1,\eta_{\xi}+2,\ldots,\eta_{\xi+1}$, there are
\[
A_{k,\kappa}=2A_{k-1,\kappa}/3=2^{\eta-\kappa+k+\sigma(\xi+1,s)}3^{\kappa-k-1}\texttt{ for }k=2,3,\ldots,\eta_{1}
\]
and for $\zeta=1,2,\ldots,\xi$ and $k=\eta_{\zeta}+1,\eta_{\zeta}+2,\ldots,\min(\eta_{\zeta+1},\kappa-1)$, if $k=\eta_{\zeta}+1$ then by the formula $3A_{k,\kappa}-2^{1+r_{\zeta}}A_{k-1,\kappa}=0$ there are
\[
A_{k,\kappa}=2^{1+r_{\zeta}}A_{k-1,\kappa}/3=2^{\eta-\kappa+k+\sigma(\xi+1,s)+\sigma(\zeta)}3^{\kappa-k-1}
\]
else by the formula $3A_{k,\kappa}-2A_{k-1,\kappa}=0$ there are
\[
A_{k,\kappa}=2A_{k-1,\kappa}/3=2^{\eta-\kappa+k+\sigma(\xi+1,s)+\sigma(\zeta)}3^{\kappa-k-1}
\]
for $k=\eta_{\zeta}+2,\eta_{\zeta}+3,\ldots,\min(\eta_{\zeta+1},\kappa-1)$. So for $\xi=1,2,\ldots,s-1$ and $\kappa=\eta_{\xi}+1,\eta_{\xi}+2,\ldots,\eta_{\xi+1}$ we obtain
\begin{equation}\label{fm:rows_l_2}
A_{k,\kappa}=2^{\eta-\kappa+k+\sigma(\xi+1,s)+\sigma(\zeta)}3^{\kappa-k-1}
\end{equation}
for $\zeta=0,1,2,\ldots,\xi$ and $k=\eta_{\zeta}+1,\eta_{\zeta}+2,\ldots,\min(\eta_{\zeta+1},\kappa-1)$.

Thus by formulas~\ref{fm:rows_l_1}-\ref{fm:rows_l_2} integers $A_{k,\kappa}$ are determined for $\kappa=2,3,\ldots,\eta$ and $k=1,2,\ldots,\kappa-1$.

Hence the matrix $\mathbf{A}^{*}$ is determined.

Let $A_{\kappa}=\sum_{j=1}^{\eta}A_{j,\kappa}$ for $\kappa=1,2,\ldots,\eta$. Then by formulas~\ref{fm:rows_r_1}-\ref{fm:rows_l_2} there are
\[
A_{\kappa}=\sum_{k=1}^{\kappa-1}2^{\eta-\kappa+k+\sigma(s)}3^{\kappa-k-1}
\]
\[
+\sum_{k=\kappa}^{\eta_{1}}2^{k-\kappa}3^{\eta-k+\kappa-1}
+\sum_{\zeta=1}^{s-1}\sum_{k=\eta_{\zeta}+1}^{\eta_{\zeta+1}}2^{k-\kappa+\sigma(\zeta)}3^{\eta-k+\kappa-1}
\]
\[
=2^{\eta-\kappa+1+\sigma(s)}3^{\kappa-1}[1-(2/3)^{\kappa-1}]+3^{\eta}[1-(2/3)^{\eta_{1}-\kappa+1}]
\]
\[
+\sum_{\zeta=1}^{s-1}2^{\eta_{\zeta}+1-\kappa+\sigma(\zeta)}3^{\eta-\eta_{\zeta}+\kappa-1}[1-(2/3)^{\eta_{\zeta+1}-\eta_{\zeta}}]
\texttt{ for }\kappa=1,2,\ldots,\eta_{1}
\]
and for $\xi=1,2,\ldots,s-1$ and $\kappa=\eta_{\xi}+1,\eta_{\xi}+2,\ldots,\min(\eta_{\xi+1},\eta)$
\[
A_{\kappa}=\sum_{\zeta=0}^{\xi}\sum_{k=\eta_{\zeta}+1}^{\min(\eta_{\zeta+1},\kappa-1)}2^{\eta-\kappa+k+\sigma(\xi+1,s)+\sigma(\zeta)}3^{\kappa-k-1}
\]
\[
+\sum_{k=\kappa}^{\eta_{\xi+1}}2^{k-\kappa}3^{\eta-k+\kappa-1}
+\sum_{\zeta=\xi+1}^{s-1}\sum_{k=\eta_{\zeta}+1}^{\eta_{\zeta+1}}
2^{k-\kappa+\sigma(\xi+1,\zeta)}3^{\eta-k+\kappa-1}
\]
\[
=\sum_{\zeta=0}^{\xi}2^{\eta-\kappa+\eta_{\zeta}+1+\sigma(\xi+1,s)+\sigma(\zeta)}3^{\kappa-\eta_{\zeta}-1}[1-(2/3)^{\min(\eta_{\zeta+1},\kappa-1)-\eta_{\zeta}}]
\]
\[
+3^{\eta}[1-(2/3)^{\eta_{\xi+1}-\kappa+1}]
+\sum_{\zeta=\xi+1}^{s-1}2^{\eta_{\zeta}+1-\kappa+\sigma(\xi+1,\zeta)}3^{\eta-\eta_{\zeta}+\kappa-1}[1-(2/3)^{\eta_{\zeta+1}-\eta_{\zeta}}].
\]

Particularly we have
\[
A_{1}=\sum_{\zeta=0}^{s-1}2^{\eta_{\zeta}+\sigma(\zeta)}3^{\eta-\eta_{\zeta}}[1-(2/3)^{\eta_{\zeta+1}-\eta_{\zeta}}].
\]
It can be rewritten as
\[
A_{1}=\sum_{\zeta=1}^{s}2^{\eta_{\zeta-1}+\sigma(\zeta-1)}3^{\eta-\eta_{\zeta}}(3^{\eta_{\zeta}-\eta_{\zeta-1}}-2^{\eta_{\zeta}-\eta_{\zeta-1}})
\]
or
\[
A_{1}=\sum_{\zeta=1}^{s-2}2^{\eta_{\zeta}+\sigma(\zeta)}3^{\eta-\eta_{\zeta}}
-\sum_{\zeta=1}^{s-2}2^{\eta_{\zeta+1}+\sigma(\zeta)}3^{\eta-\eta_{\zeta+1}}
\]
\[
+3^{\eta}-2^{\eta_{1}}3^{\eta-\eta_{1}}
+2^{\eta_{s-1}+\sigma(s-1)}3^{\eta_{s}-\eta_{s-1}}-2^{\eta+\sigma(s-1)}
\]
\[
=\sum_{\zeta=1}^{s}2^{\eta_{\zeta}+\sigma(\zeta-1)}3^{\eta-\eta_{\zeta}}(2^{\sigma(\zeta,\zeta)}-1)+Det\mathbf{A}.
\]

This completes the proof of the lemma.
\end{proof}

\begin{lemma}\label{le:diophantine_eq_sol}
Given a number
\[
A_{1}=\sum_{\zeta=1}^{s}2^{\eta_{\zeta}+\sigma(\zeta-1)}3^{\eta-\eta_{\zeta}}(2^{\sigma(\zeta,\zeta)}-1)+Det\mathbf{A}
\]
where $Det\mathbf{A}=3^{\eta}-2^{\eta+\sigma(s)}$, and for $\xi=1,2,\cdots,s$ there are $\eta=\eta_{s}$, $\eta_{0}=0$,
\[
\sigma(\zeta)=\sum_{j=1}^{\zeta}r_{j}\texttt{, }\sigma(\xi,\zeta)=\sum_{j=\xi}^{\zeta}r_{j}\texttt{, and }\eta_{\xi}=\sum_{j=0}^{\xi-1}(m_{j}-i_{j}+1)\texttt{ where }0\leq i_{j}\leq m_{j},
\]
$A_{1}$ can not be exactly dividable by $-Det\mathbf{A}$ for $1<s<\infty$.
\end{lemma}

\begin{proof}
Let rewrite $A_{1}=2^{\eta_{1}}A_{1}^{*}+Det\mathbf{A}$ where by the assumption $r_{c}(m,t,i,s)=1$ there must be $\eta>s$, $A_{1}$ and also $A_{1}^{*}$ must be exactly dividable by $-Det\mathbf{A}$, and
\begin{equation}\label{fm:rows_kappa_1_d}
A_{1}^{*}=\sum_{\zeta=1}^{s}2^{\eta_{\zeta}-\eta_{1}+\sigma(\zeta-1)}3^{\eta-\eta_{\zeta}}(2^{\sigma(\zeta,\zeta)}-1)=-Det\mathbf{A}\times n^{*}
\end{equation}
where $n^{*}$ is a positive odd number.

Let $a_{\zeta}^{*}=2^{\eta_{\zeta}-\eta_{1}+\sigma(\zeta-1)}3^{\eta-\eta_{\zeta}}$ and consider the Diophantine equation
\begin{equation}\label{eq:diophantine_s}
\sum_{\zeta=1}^{s}a_{\zeta}^{*}y_{\zeta}^{*}=-Det\mathbf{A}\times n^{*}\texttt{ where }(a_{1}^{*},a_{2}^{*},\ldots,a_{s}^{*})=1.
\end{equation}

When $s=2$, by Lemma~\ref{le:rows_a_1_s2} there must be
\[
r_{1}=r_{2}\texttt{, }2^{\eta_{1}+\sigma(s,s)}=3^{\eta-\eta_{1}}+1\texttt{, and }2^{\eta-\eta_{1}+\sigma(1)}=3^{\eta_{1}}+1.
\]
Thus there shall be $\eta-\eta_{1}=\eta_{1}=1$ and $\sigma(j,j)=r_{j}=1\texttt{ for }j=1,2$, so that there are  $\eta=s$ and $n_{i_{j}}(m_{j},t_{j})=n_{0}(0,1)=1$ for $j=0,1,2$, i.e., no nontrivial Collatz cycle exists in $S_{c}(s)$ and there shall be $s=1$. It is a contradiction. Thus $A_{1}^{*}$ and also $A_{1}$ can not be exactly dividable by $-Det\mathbf{A}$ when $s=2$.

When $s=3$, by Lemma~\ref{le:rows_a_1_s3} there must be
\[
r_{1}=r_{2}=r_{3}\texttt{, }2^{\eta_{1}+\sigma(s,s)}=3^{\eta-\eta_{2}}+1\texttt{, }
\]
\[
2^{\eta-\eta_{2}+\sigma(2,2)}=3^{\eta_{2}-\eta_{1}}+1\texttt{, and }2^{\eta_{2}-\eta_{1}+\sigma(1)}=3^{\eta_{1}}+1.
\]
Thus there shall be $\eta-\eta_{2}=\eta_{2}-\eta_{1}=\eta_{1}=1$ and $\sigma(j,j)=r_{j}=1\texttt{ for }j=1,2,3$, so that there are $\eta=s$ and $n_{i_{j}}(m_{j},t_{j})=n_{0}(0,1)=1$ for $j=0,1,2,3$, i.e., no nontrivial Collatz cycle exists in $S_{c}(s)$ and there shall be $s=1$. It is a contradiction. Thus $A_{1}^{*}$ and also $A_{1}$ can not be exactly dividable by $-Det\mathbf{A}$ when $s=3$.

When $3<s<\infty$, let $\eta_{\zeta}-\eta_{\zeta-1}=1$ and $\sigma(\zeta,\zeta)=r_{\zeta}=1$ for $\zeta=1,2,\ldots,s$. Then a group of solutions of Diophantine Equation~\ref{eq:diophantine_s} are $y_{\zeta}^{*}=2^{\sigma(\zeta,\zeta)}-1=n^{*}=1$ for $\zeta=1,2,\ldots,s$ and there are
\begin{equation}\label{fm:rows_a_1_s}
r_{1}=r_{2}=\cdots=r_{s}\texttt{, }2^{\eta_{1}+\sigma(s,s)}=3^{\eta-\eta_{s-1}}+1\texttt{, and}
\end{equation}
\[
2^{\eta_{\zeta+1}-\eta_{\zeta}+\sigma(\zeta,\zeta)}=3^{\eta_{\zeta}-\eta_{\zeta-1}}+1\texttt{ for }\zeta=s-1,s-2,\ldots,2,1.
\]
Because the group of values $y_{\zeta}^{*}=2^{\sigma(\zeta,\zeta)}-1=n^{*}=1$ for $\zeta=1,2,\ldots,s$ are a group of particular solutions of Diophantine Equation~\ref{eq:diophantine_s}, by the theory of Diophantine Equations~\cite{Z}, any other group of solutions of Diophantine Equation~\ref{eq:diophantine_s} do not hold since at least one $y_{\zeta}^{*}$ in the other group of solutions is less than or equal to zero where $1\leq\zeta\leq s$, which conflicts with the fact that there should be $y_{\zeta}^{*}=2^{\sigma(\zeta,\zeta)}-1\geq1$ for $\zeta=1,2,\ldots,s$ in any group of solutions of Diophantine Equation~\ref{eq:diophantine_s}, so that there must be $2^{\sigma(\zeta,\zeta)}-1=n^{*}$ for $\zeta=1,2,\ldots,s$ and Formula~\ref{fm:rows_a_1_s} must hold.
Thus there must be $\eta_{\zeta}-\eta_{\zeta-1}=1$ and $\sigma(\zeta,\zeta)=r_{\zeta}=1\texttt{ for }\zeta=1,2,\ldots,s$, so that there are $\eta=s$ and $n_{i_{j}}(m_{j},t_{j})=n_{0}(0,1)=1$ for $j=0,1,2,\ldots,s$, i.e., no nontrivial Collatz cycle exists in $S_{c}(s)$ and there shall be $s=1$. It is a contradiction. Thus $A_{1}^{*}$ and also $A_{1}$ can not be exactly dividable by $-Det\mathbf{A}$ when $3<s<\infty$.

Hence the number $A_{1}$ can not be exactly dividable by $-Det\mathbf{A}$ for $1<s<\infty$.

This completes the proof of the lemma.
\end{proof}

\begin{lemma}\label{le:matrix_eq_sol}
Given a number
\begin{equation}\label{fm:rows_a_1}
A_{1}=\sum_{\zeta=1}^{s}2^{\eta_{\zeta-1}+\sigma(\zeta-1)}3^{\eta-\eta_{\zeta}}(3^{\eta_{\zeta}-\eta_{\zeta-1}}-2^{\eta_{\zeta}-\eta_{\zeta-1}})
\end{equation}
\[
=\sum_{\zeta=1}^{s}2^{\eta_{\zeta}+\sigma(\zeta-1)}3^{\eta-\eta_{\zeta}}(2^{\sigma(\zeta,\zeta)}-1)+Det\mathbf{A}
\]
where $Det\mathbf{A}=3^{\eta}-2^{\eta+\sigma(s)}$, and for $\xi=1,2,\cdots,s$ there are $\eta=\eta_{s}$, $\eta_{0}=0$,
\[
\sigma(\zeta)=\sum_{j=1}^{\zeta}r_{j}\texttt{, }\sigma(\xi,\zeta)=\sum_{j=\xi}^{\zeta}r_{j}\texttt{, and }\eta_{\xi}=\sum_{j=0}^{\xi-1}(m_{j}-i_{j}+1)\texttt{ where }0\leq i_{j}\leq m_{j}
\]
and $\eta_{j+1}-\eta_{j}=m_{j}-i_{j}+1$, $A_{1}$ can not be exactly dividable by $-Det\mathbf{A}$ for $s>1$.
\end{lemma}

\begin{proof}
Let rewrite $A_{1}=2^{\eta_{1}}A_{1}^{*}+Det\mathbf{A}$ where by the assumption $r_{c}(m,t,i,s)=1$ there must be $\eta>s$, $A_{1}$ and also $A_{1}^{*}$ must be exactly dividable by $-Det\mathbf{A}$, and
\begin{equation}\label{fm:rows_kappa_1}
A_{1}^{*}=\sum_{\zeta=1}^{s}2^{\eta_{\zeta}-\eta_{1}+\sigma(\zeta-1)}3^{\eta-\eta_{\zeta}}(2^{\sigma(\zeta,\zeta)}-1).
\end{equation}

Let $a_{\zeta}=2^{\eta_{\zeta-1}+\sigma(\zeta-1)}3^{\eta-\eta_{\zeta}}$ and $a_{\zeta}^{*}=2^{\eta_{\zeta}-\eta_{1}+\sigma(\zeta-1)}3^{\eta-\eta_{\zeta}}$ for $\zeta=1,2,\ldots,s$. Then we consider equations
\begin{equation}\label{eq:rows_a_1}
\sum_{\zeta=1}^{s}a_{\zeta}y_{\zeta}=-Det\mathbf{A}\times n\texttt{ where }(a_{1},a_{2},\ldots,a_{s})=1
\end{equation}
and
\begin{equation}\label{eq:rows_kappa_1}
\sum_{\zeta=1}^{s}a_{\zeta}^{*}y_{\zeta}^{*}=-Det\mathbf{A}\times n^{*}\texttt{ where }(a_{1}^{*},a_{2}^{*},\ldots,a_{s}^{*})=1.
\end{equation}

If $A_{1}=-Det\mathbf{A}\times n$ and $A_{1}^{*}=-Det\mathbf{A}\times n^{*}$, then there shall be $n=2^{\eta_{1}}n^{*}-1$ where $n$ and $n^{*}$ must be positive odd numbers.

If $A_{1}=-Det\mathbf{A}\times n$, then a group of positive odd numbers
\[
y_{\zeta}=3^{\eta_{\zeta}-\eta_{\zeta-1}}-2^{\eta_{\zeta}-\eta_{\zeta-1}}=3^{m_{\zeta-1}-i_{\zeta-1}+1}-2^{m_{\zeta-1}-i_{\zeta-1}+1}\texttt{ for }\zeta=1,2,\ldots,s
\]
must be a solution of Equation~\ref{eq:rows_a_1}.

Because each $a_{\zeta}$ is a product $\alpha_{\zeta}\beta_{\zeta}$ where
\[
\alpha_{\zeta}=2^{\eta_{\zeta-1}+\sigma(\zeta-1)}\texttt{ and }\beta_{\zeta}=3^{\eta-\eta_{\zeta}}\texttt{ for }\zeta=1,2,\ldots,s,
\]
any two different products $a_{i}$ and $a_{j}$ can not be summed as a product $a_{k}$, i.e.,
\[
a_{i}+a_{j}\neq a_{k}\texttt{ or }\alpha_{i}\beta_{i}+\alpha_{j}\beta_{j}\neq\alpha_{k}\beta_{k}\texttt{ where }1\leq i,j,k\leq s\texttt{ and }i\neq j.
\]
Thus there must be $y_{i}\neq y_{j}$ when $1\leq i<j\leq s$.

If $A_{1}^{*}=-Det\mathbf{A}\times n^{*}$, then a group of positive odd numbers $y_{\zeta}^{*}=2^{\sigma(\zeta,\zeta)}-1$ for $\zeta=1,2,\ldots,s$ must be a solution of Equation~\ref{eq:rows_kappa_1}.


Because each $a_{\zeta}^{*}$ is a product $\alpha_{\zeta}^{*}\beta_{\zeta}$ where
\[
\alpha_{\zeta}^{*}=2^{\eta_{\zeta}-\eta_{1}+\sigma(\zeta-1)}\texttt{ and }\beta_{\zeta}=3^{\eta-\eta_{\zeta}}\texttt{ for }\zeta=1,2,\ldots,s,
\]
any two different products $a_{i}^{*}$ and $a_{j}^{*}$ can not be summed as a product $a_{k}^{*}$, i.e.,
\[
a_{i}^{*}+a_{j}^{*}\neq a_{k}^{*}\texttt{ or }\alpha_{i}^{*}\beta_{i}+\alpha_{j}^{*}\beta_{j}\neq\alpha_{k}^{*}\beta_{k}\texttt{ where }1\leq i,j,k\leq s\texttt{ and }i\neq j.
\]
Thus there must be $y_{i}^{*}\neq y_{j}^{*}$ when $1\leq i<j\leq s$.

Then we shall prove that no group of numbers
\[
y_{\zeta}=3^{m_{\zeta-1}-i_{\zeta-1}+1}-2^{m_{\zeta-1}-i_{\zeta-1}+1}\texttt{ for }\zeta=1,2,\ldots,s
\]
are the solution of Equation~\ref{eq:rows_a_1} and also no group of numbers
\[
y_{\zeta}^{*}=2^{\sigma(\zeta,\zeta)}-1\texttt{ for }\zeta=1,2,\ldots,s
\]
are the solution of Equation~\ref{eq:rows_kappa_1} when $s>2$.

First, when at least one pair of nonnegative integers $m_{\xi-1}-i_{\xi-1}$ and $m_{\zeta-1}-i_{\zeta-1}$ satisfy $m_{\xi-1}-i_{\xi-1}=m_{\zeta-1}-i_{\zeta-1}$ where $1\leq\xi<\zeta\leq s$, no group of numbers
\[
y_{\zeta}=3^{m_{\zeta-1}-i_{\zeta-1}+1}-2^{m_{\zeta-1}-i_{\zeta-1}+1}\texttt{ for }\zeta=1,2,\ldots,s
\]
are the solution of Equation~\ref{eq:rows_a_1} so that $A_{1}\neq-Det\mathbf{A}\times n$.

By Lemma~\ref{le:odd_seqs_m} one half of nonnegative integers $m_{\zeta-1}-i_{\zeta-1}$ satisfy $m_{\zeta-1}-i_{\zeta-1}=0$ and one quarter of nonnegative integers $m_{\zeta-1}-i_{\zeta-1}$ satisfy $m_{\zeta-1}-i_{\zeta-1}=1$ and so on for $\zeta=1,2,\ldots,s$. Then more than $s/2-1$ pairs of nonnegative integers $m_{\xi-1}-i_{\xi-1}$ and $m_{\zeta-1}-i_{\zeta-1}$ satisfy $m_{\xi-1}-i_{\xi-1}=m_{\zeta-1}-i_{\zeta-1}$ where $1\leq\xi<\zeta\leq s$. Thus $A_{1}$ can not be exactly dividable by $-Det\mathbf{A}$ when $s>3$.

Second, when at least one pair of positive integers $r_{i}$ and $r_{j}$ satisfy $r_{i}=r_{j}$ where $1\leq i<j\leq s$, no group of numbers $y_{\zeta}^{*}=2^{\sigma(\zeta,\zeta)}-1\texttt{ for }\zeta=1,2,\ldots,s$ are the solution of Equation~\ref{eq:rows_kappa_1} so that $A_{1}^{*}\neq-Det\mathbf{A}\times n^{*}$.

By Lemma~\ref{le:odd_seq_r1} one half of positive integers $\sigma(\zeta,\zeta)=r_{\zeta}$ satisfy $r_{\zeta}=1$ and one quarter of positive integers $\sigma(\zeta,\zeta)=r_{\zeta}$ satisfy $r_{\zeta}=2$ and so on for $\zeta=1,2,\ldots,s$. Then more than $s/2-1$ pairs of positive integers $r_{i}$ and $r_{j}$ satisfy $r_{i}=r_{j}$ where $1\leq i<j\leq s$. Thus $A_{1}^{*}$ and also $A_{1}$ can not be exactly dividable by $-Det\mathbf{A}$ when $s>3$.

Third, by combining these two situations above, at least $s-2$ equalities $y_{i}=y_{j}$ and/or $y_{i}^{*}=y_{j}^{*}$ hold where $1\leq i<j\leq s$. Thus $A_{1}$ can not be exactly dividable by $-Det\mathbf{A}$ when $s>2$.

On the other hand, when $s=3$, Formula~\ref{fm:rows_kappa_1} becomes
\[
A_{1}^{*}=3^{\eta-\eta_{1}}(2^{\sigma(1,1)}-1)
+2^{\eta_{2}-\eta_{1}+\sigma(1)}3^{\eta-\eta_{2}}(2^{\sigma(2,2)}-1)
+2^{\eta-\eta_{1}+\sigma(2)}(2^{\sigma(3,3)}-1).
\]
By Lemma~\ref{le:rows_a_1_s3} if $A_{1}^{*}$ can be exactly dividable by $-Det\mathbf{A}$ then there must be $r_{1}=r_{2}=r_{3}$. But it contradicts with the fact that by the analysis above, if one pair of positive integers $r_{i}$ and $r_{j}$ satisfy $r_{i}=r_{j}$ where $1\leq i<j\leq s$ then no group of numbers $y_{\zeta}^{*}=2^{\sigma(\zeta,\zeta)}-1$ for $\zeta=1,2,3$ are the solution of Equation~\ref{eq:rows_kappa_1} so that $A_{1}^{*}\neq-Det\mathbf{A}\times n^{*}$. Furthermore, when $r_{1}=r_{2}=r_{3}$, Formula~\ref{fm:rows_a_1_s3} and Formula~\ref{fm:rows_a_1_v3} hold, so that there are $\eta=s$ and $n_{i_{j}}(m_{j},t_{j})=n_{0}(0,1)=1$ for $j=0,1,2,3$, i.e., no nontrivial Collatz cycle exists in $S_{c}(s)$ and there shall be $s=1$. Thus $A_{1}^{*}$ and also $A_{1}$ can not be exactly dividable by $-Det\mathbf{A}$ when $s=3$.

Finally, when $s=2$, Formula~\ref{fm:rows_kappa_1} becomes
\[
A_{1}^{*}=3^{\eta-\eta_{1}}(2^{\sigma(1,1)}-1)+2^{\eta-\eta_{1}+\sigma(1)}(2^{\sigma(2,2)}-1).
\]
By Lemma~\ref{le:rows_a_1_s2} if $A_{1}^{*}$ can be exactly dividable by $-Det\mathbf{A}$ then there must be $r_{1}=r_{2}$. But it contradicts with the fact that by the analysis above, if one pair of positive integers $r_{1}$ and $r_{2}$ satisfy $r_{1}=r_{2}$ then no group of numbers $y_{\zeta}^{*}=2^{\sigma(\zeta,\zeta)}-1$ for $\zeta=1,2$ are the solution of Equation~\ref{eq:rows_kappa_1} so that $A_{1}^{*}\neq-Det\mathbf{A}\times n^{*}$. Furthermore, when $r_{1}=r_{2}$, Formula~\ref{fm:rows_a_1_s2} holds and there shall be $\eta-\eta_{1}=\eta_{1}=1$, $r_{1}=\sigma(1,1)=1$, and $r_{2}=\sigma(s,s)=1$, so that there are $\eta=s$ and $n_{i_{j}}(m_{j},t_{j})=n_{0}(0,1)=1$ for $j=0,1,2$, i.e., no nontrivial Collatz cycle exists in $S_{c}(s)$ and there shall be $s=1$. Thus $A_{1}^{*}$ and also $A_{1}$ can not be exactly dividable by $-Det\mathbf{A}$ when $s=2$.

Hence the number $A_{1}$ can not be exactly dividable by $-Det\mathbf{A}$ for $1<s<\infty$.

This completes the proof of the lemma.
\end{proof}

\begin{lemma}\label{le:odd_seqs_m}
Given a series of Collatz's Sequences
\[
S_{c}(s)=\{S_{c}(m_{j},t_{j})|j=0,1,2,\ldots,s-1\}
\]
where $m_{0}=m$, $t_{0}=t$, and for $j=0,1,2,\ldots,s-1$ the odd number $n_{i_{j}}(m_{j},t_{j})$ is greater than 1 and the even number
\[
n_{m_{j}+1}(m_{j},t_{j})=2^{r_{j+1}}n_{i_{j+1}}(m_{j+1},t_{j+1})\texttt{ where }r_{j+1}\geq1\texttt{ and }0\leq i_{j+1}\leq m_{j+1}.
\]

Then $1/2^{k}$ of positive integers $m_{j}-i_{j}+1$ satisfy $m_{j}-i_{j}+1=k$ for $k=1,2,\ldots$ as $s$ tends to infinity.
\end{lemma}

\begin{proof}
Let Collatz's Sequences $S_{c}(m_{j},t_{j})$ for $j=0,1,2,\ldots,s-1$ be samples of the Collatz's Sequence $S_{c}(m,t)$. Then $s$ is the number of independent samples of the Collatz's Sequence $S_{c}(m,t)$ in the series of Collatz's Sequences $S_{c}(s)$.

By Lemma~\ref{le:odd_seq} there are $s$ even numbers
\[
n_{i_{j}}(m_{j},t_{j})+1=2^{m_{j}+1-i_{j}}3^{i_{j}}t_{j}\texttt{ for }j=0,1,2,\ldots,s-1.
\]

By Lemma~\ref{le:odd_seq_m}, there are $m_{j}-i_{j}=0$ for $n_{i_{j}}(m_{j},t_{j})=4k+1$ where $i_{j}=m_{j}=0$ or $i_{j}=m_{j}$, and $m_{j}-i_{j}>0$ for $n_{i_{j}}(m_{j},t_{j})=4k+3$ where $0\leq i_{j}<m_{j}$. Then there are $s/2$ odd numbers $n_{i_{j}}(m_{j},t_{j})$ with $m_{j}-i_{j}=0$.

For another $s/2$ odd numbers $n_{i_{j}}(m_{j},t_{j})=4k+3$, because the even number
\[
n_{i_{j}}(m_{j},t_{j})+1=2^{m_{j}+1-i_{j}}3^{i_{j}}t_{j}=4(k+1),
\]
there are $2^{m_{j}-1-i_{j}}3^{i_{j}}t_{j}=k+1$ and
\[
m_{j}-i_{j}=\kappa+1\texttt{ where }k+1=2^{\kappa}3^{i_{j}}t_{j}\texttt{ for }\kappa=0,1,2,\ldots.
\]
Then there are $s/2^{2}$ odd numbers $n_{0}(m_{j},t_{j})$ with $m_{j}-i_{j}=1$, $s/2^{3}$ odd numbers $n_{0}(m_{j},t_{j})$ with $m_{j}-i_{j}=2$, and $s/2^{4}$ odd numbers $n_{0}(m_{j},t_{j})$ with $m_{j}-i_{j}=3,\ldots$.

Therefore as $s$ tends to infinity, i.e., the number of independent samples of the Collatz's Sequence $S_{c}(m,t)$ tends to infinity, by the law of great numbers, the distribution of the $s$ even numbers $n_{i_{j}}(m_{j},t_{j})+1$ is similar to the distribution of even numbers $n_{e}(\gamma,\lambda)$ stated by Lemma~\ref{le:even_distribution}, so that $1/2^{k}$ of positive integers $m_{j}-i_{j}+1$ satisfy $m_{j}-i_{j}+1=k$ for $k=1,2,\ldots$ as $s$ tends to infinity.

This completes the proof of the lemma.
\end{proof}

\begin{lemma}\label{le:rows_a_1_s2}
For $s=2$, given a number
\[
A_{1}^{*}=3^{\eta-\eta_{1}}(2^{\sigma(1,1)}-1)+2^{\eta-\eta_{1}+\sigma(1)}(2^{\sigma(2,2)}-1)
\]
with Definition~\ref{def:sigma_eta}, if and only if
\begin{equation}\label{fm:rows_a_1_s2}
r_{1}=r_{2}\texttt{, }2^{\eta_{1}+\sigma(s,s)}=3^{\eta-\eta_{1}}+1\texttt{, and }2^{\eta-\eta_{1}+\sigma(1)}=3^{\eta_{1}}+1,
\end{equation}
then $A_{1}^{*}$ can be exactly dividable by $-Det\mathbf{A}$ where $Det\mathbf{A}=3^{\eta}-2^{\eta+\sigma(s)}$.
\end{lemma}

\begin{proof}
If Formula~\ref{fm:rows_a_1_s2} holds, then there shall be $\eta-\eta_{1}=\eta_{1}=1$, $\sigma(s,s)=\sigma(1)=1$, and $A_{1}^{*}=-Det\mathbf{A}(2^{\sigma(1,1)}-1)$, so that $A_{1}^{*}$ can be exactly dividable by $-Det\mathbf{A}$.

On the other hand, if $A_{1}^{*}=-Det\mathbf{A}\times n^{*}$ where $n^{*}$ is odd, then there shall be
\[
3^{\eta-\eta_{1}}(2^{\sigma(1,1)}-1+3^{\eta_{1}}n^{*})=2^{\eta-\eta_{1}+\sigma(1)}(2^{\eta_{1}+\sigma(s,s)}n^{*}-2^{\sigma(2,2)}+1).
\]
Let $a_{1}^{*}=3^{\eta-\eta_{1}}$ and $a_{2}^{*}=2^{\eta-\eta_{1}+\sigma(1)}$. Then we consider the Diophantine equation
\begin{equation}\label{eq:rows_a_1_s2}
a_{1}^{*}y_{1}^{*}+a_{2}^{*}y_{2}^{*}=-Det\mathbf{A}\times n^{*}\texttt{ where }(a_{1}^{*},a_{2}^{*})=1.
\end{equation}
Because the group of values $y_{\zeta}^{*}=2^{\sigma(\zeta,\zeta)}-1=n^{*}=1$ for $\zeta=1,2$ are a group of particular solutions of Diophantine Equation~\ref{eq:rows_a_1_s2}, by the theory of Diophantine Equations~\cite{Z}, any other group of solutions of Diophantine Equation~\ref{eq:rows_a_1_s2} do not hold since at least one $y_{\zeta}^{*}$ in the other group of solutions is less than or equal to zero where $1\leq\zeta\leq2$, which conflicts with the fact that there should be $y_{\zeta}^{*}=2^{\sigma(\zeta,\zeta)}-1\geq1$ for $\zeta=1,2$ in any group of solutions of Diophantine Equation~\ref{eq:rows_a_1_s2}, so that there must be $2^{\sigma(\zeta,\zeta)}-1=n^{*}$ for $\zeta=1,2$ and Formula~\ref{fm:rows_a_1_s2} must hold. Then there shall be $n^{*}=1$, $\eta-\eta_{1}=\eta_{1}=1$, and $r_{1}=r_{2}=\sigma(1,1)=\sigma(2,2)=1$.

So if and only if Formula~\ref{fm:rows_a_1_s2} holds then $A_{1}^{*}$ can be exactly dividable by $-Det\mathbf{A}$.

This completes the proof of the lemma.
\end{proof}

\begin{lemma}\label{le:rows_a_1_s3}
For $s=3$, given a number
\[
A_{1}^{*}=3^{\eta-\eta_{1}}(2^{\sigma(1,1)}-1)
+2^{\eta_{2}-\eta_{1}+\sigma(1)}3^{\eta-\eta_{2}}(2^{\sigma(2,2)}-1)
\]
\[
+2^{\eta-\eta_{1}+\sigma(2)}(2^{\sigma(3,3)}-1)
\]
with Definition~\ref{def:sigma_eta}, if and only if
\begin{equation}\label{fm:rows_a_1_s3}
r_{1}=r_{2}=r_{3}\texttt{, }2^{\eta_{1}+\sigma(s,s)}=3^{\eta-\eta_{2}}+1\texttt{, }
\end{equation}
\[
2^{\eta-\eta_{2}+\sigma(2,2)}=3^{\eta_{2}-\eta_{1}}+1\texttt{, and }2^{\eta_{2}-\eta_{1}+\sigma(1)}=3^{\eta_{1}}+1,
\]
then $A_{1}^{*}$ can be exactly dividable by $-Det\mathbf{A}$ where $Det\mathbf{A}=3^{\eta}-2^{\eta+\sigma(s)}$.
\end{lemma}

\begin{proof}
If Formula~\ref{fm:rows_a_1_s3} holds, then there shall be $A_{1}^{*}=-Det\mathbf{A}(2^{\sigma(1,1)}-1)$ where
\begin{equation}\label{fm:rows_a_1_v3}
\eta-\eta_{2}=\eta_{2}-\eta_{1}=\eta_{1}=1\texttt{ and }\sigma(j,j)=r_{j}=1\texttt{ for }j=1,2,3.
\end{equation}

On the other hand, if $A_{1}^{*}=-Det\mathbf{A}\times n^{*}$ where $n^{*}$ is odd, then there shall be
\[
 3^{\eta-\eta_{2}}[3^{\eta_{2}-\eta_{1}}(2^{\sigma(1,1)}-1+3^{\eta_{1}}n^{*})
+2^{\eta_{2}-\eta_{1}+\sigma(1)}(2^{\sigma(2,2)}-1)]
\]
\[
=2^{\eta-\eta_{1}+\sigma(2)}(2^{\eta_{1}+\sigma(s,s)}n^{*}-2^{\sigma(3,3)}+1).
\]
Let $a_{\zeta}^{*}=2^{\eta_{\zeta}-\eta_{1}+\sigma(\zeta-1)}3^{\eta-\eta_{\zeta}}$ for $\zeta=1,2,3$ and consider the Diophantine equation
\begin{equation}\label{eq:rows_a_1_s3}
a_{1}^{*}y_{1}^{*}+a_{2}^{*}y_{2}^{*}+a_{3}^{*}y_{3}^{*}=-Det\mathbf{A}\times n^{*}\texttt{ where }(a_{1}^{*},a_{2}^{*},a_{3}^{*})=1.
\end{equation}
Because the group of values $y_{\zeta}^{*}=2^{\sigma(\zeta,\zeta)}-1=n^{*}=1$ for $\zeta=1,2,3$ are a group of particular solutions of Diophantine Equation~\ref{eq:rows_a_1_s3}, by the theory of Diophantine Equations~\cite{Z}, any other group of solutions of Diophantine Equation~\ref{eq:rows_a_1_s3} do not hold since at least one $y_{\zeta}^{*}$ in the other group of solutions is less than or equal to zero where $1\leq\zeta\leq3$, which conflicts with the fact that there should be $y_{\zeta}^{*}=2^{\sigma(\zeta,\zeta)}-1\geq1$ for $\zeta=1,2,3$ in any group of solutions of Diophantine Equation~\ref{eq:rows_a_1_s3}, so that there must be $2^{\sigma(\zeta,\zeta)}-1=n^{*}$ for $\zeta=1,2,3$ and Formula~\ref{fm:rows_a_1_s3} must hold. Then there shall be $n^{*}=1$ and Formula~\ref{fm:rows_a_1_v3} holds.

So if and only if Formula~\ref{fm:rows_a_1_s3} holds then $A_{1}^{*}$ can be exactly dividable by $-Det\mathbf{A}$.

This completes the proof of the lemma.
\end{proof}

\begin{lemma}[Distribution of even numbers]\label{le:even_distribution}
As a natural number $s_{e}$ tends to infinite, there are $s_{e}/2^{\gamma}$ even numbers $n_{e}(\gamma,\lambda)\leq2s_{e}$ for $\gamma=1,2,\ldots$ where
\[
n_{e}(\gamma,\lambda)=2^{\gamma}(2\lambda+1)\texttt{ where }\gamma=1,2,\ldots\texttt{ and }\lambda=0,1,2,\ldots.
\]
\end{lemma}

\begin{proof}
Since all even numbers can be expressed as
\[
n_{e}(\gamma,\lambda)=2^{\gamma}(2\lambda+1)\texttt{ where }\gamma=1,2,\ldots\texttt{ and }\lambda=0,1,2,\ldots,
\]
then as a natural number $s_{e}$ tends to infinite, there are $s_{e}/2^{\gamma}$ even numbers $n_{e}(\gamma,\lambda)\leq2s_{e}$ for $\gamma=1,2,\ldots$.

This completes the proof of the lemma.
\end{proof}

\begin{lemma}\label{le:odd_seqs}
Given a nonnegative integer $m$ and an odd number $t$ not divisible by $3$ with $m_{0}=m$ and $t_{0}=t$, there exist a unique series of Collatz's Sequences
\[
S_{c}(s)=\{S_{c}(m_{j},t_{j})|j=0,1,2,\ldots,s-1\}
\]
where $n_{i_{j}}(m_{j},t_{j})$ is odd and greater than 1, $0\leq i_{j}\leq m_{j}$, and the even number
\[
n_{m_{j}+1}(m_{j},t_{j})=2^{r_{j+1}}n_{i_{j+1}}(m_{j+1},t_{j+1})\texttt{ where }r_{j+1}\geq1\texttt{ and }0\leq i_{j+1}\leq m_{j+1},
\]
for $j=0,1,2,\ldots,s-1$, and the odd number $n_{i_{s}}(m_{s},t_{s})$ is equal to 1 for $s$ finite.
\end{lemma}

\begin{proof}
First, by Lemma~\ref{le:odd_seq}, given a nonnegative integer $m$ and an odd number $t$ not divisible by $3$, there exists a unique finite Collatz's Sequence
\[
S_{c}(m,t)=\{n_{0}(m,t),n_{1}(m,t),n_{2}(m,t),\ldots,n_{m}(m,t),n_{m+1}(m,t)\}
\]
where $n_{0}(m,t)$ is odd, each $n_{i}(m,t)=(3n_{i-1}(m,t)+1)/2$ for $i=1,2,\ldots,m$ is odd, and the natural number $n_{m+1}(m,t)=(3n_{m}(m,t)+1)/2$ is even.

Second, let $m_{0}=m$, $t_{0}=t$, and the even number
\[
n_{m+1}(m,t)=2^{r_{1}}n_{i_{1}}(m_{1},t_{1})\texttt{ where }r_{1}\geq1\texttt{ and }0\leq i_{1}\leq m_{1}.
\]
Then when $n_{i_{1}}(m_{1},t_{1})$ is equal to 1, we form a series of Collatz's Sequences
\[
S_{c}(s)=\{S_{c}(m_{j},t_{j})|j=0,1,2,\ldots,s-1\}\texttt{ where }s=1,
\]
otherwise since by Lemma~\ref{le:odd_seq_1_c} no nontrivial Collatz cycle exists in any Collatz's Sequence $S_{c}(m,t)$, by Lemma~\ref{le:odd_seq} and given a pair of the nonnegative integer $m_{1}$ and the odd number $t_{1}$, there exist a unique finite Collatz's Sequence $S_{c}(m_{1},t_{1})$ and an even number
\[
n_{m_{1}+1}(m_{1},t_{1})=2^{r_{2}}n_{i_{2}}(m_{2},t_{2})\texttt{ where }r_{2}\geq1\texttt{ and }0\leq i_{2}\leq m_{2}.
\]
If $n_{i_{2}}(m_{2},t_{2})$ is equal to 1 then we form a series of Collatz's Sequences
\[
S_{c}(s)=\{S_{c}(m_{j},t_{j})|j=0,1,2,\ldots,s-1\}\texttt{ where }s=2.
\]

Third, similar to the above way, for $j=2,3,\ldots$, if $n_{i_{j}}(m_{j},t_{j})$ is greater than 1 where $0\leq i_{j}\leq m_{j}$, then since by Lemma~\ref{le:odd_seq_s_c} no nontrivial Collatz cycle exists in any series of Collatz's Sequences $S_{c}(s)$ where $0<s<\infty$, by Lemma~\ref{le:odd_seq}, given a pair of the nonnegative integer $m_{j}$ and the odd number $t_{j}$, there exist a unique finite Collatz's Sequence $S_{c}(m_{j},t_{j})$ and an even number
\[
n_{m_{j}+1}(m_{j},t_{j})=2^{r_{j+1}}n_{i_{j+1}}(m_{j+1},t_{j+1})\texttt{ where }r_{j+1}\geq1\texttt{ and }0\leq i_{j+1}\leq m_{j+1}.
\]
When $n_{i_{j+1}}(m_{j+1},t_{j+1})$ is equal to 1, we form a series of Collatz's Sequences
\[
S_{c}(s)=\{S_{c}(m_{i},t_{i})|i=0,1,2,\ldots,s-1\}\texttt{ where }s=j+1.
\]

Finally, we obtain a unique series of Collatz's Sequences
\[
S_{c}(s)=\{S_{c}(m_{j},t_{j})|j=0,1,2,\ldots,s-1\}\texttt{ where }m_{0}=m\texttt{ and }t_{0}=t,
\]
the odd number $n_{i_{j}}(m_{j},t_{j})$ is greater than 1 for $j=0,1,2,\ldots,s-1$, and the even number
\[
n_{m_{j}+1}(m_{j},t_{j})=2^{r_{j+1}}n_{i_{j+1}}(m_{j+1},t_{j+1})\texttt{ where }r_{j+1}\geq1\texttt{ and }0\leq i_{j+1}\leq m_{j+1}.
\]
If $s$ is finite or $s<\infty$ then the odd number $n_{i_{s}}(m_{s},t_{s})$ is equal to 1.

This completes the proof of the lemma.
\end{proof}

\begin{lemma}\label{le:odd_seq_s_e}
No nontrivial Collatz cycle exists in a series of Collatz's Sequences
\[
S_{c}(s)=\{S_{c}(m_{j},t_{j})|j=0,1,2,\ldots,s-1\}
\]
and the odd number $n_{i_{s}}(m_{s},t_{s})$ in the series of Collatz's Sequences $S_{c}(s)$ is equal to 1 as $s$ tends to infinite where $0\leq i_{s}\leq m_{s}$.
\end{lemma}

\begin{proof}
Given a series of Collatz's Sequences
\[
S_{c}(s)=\{S_{c}(m_{j},t_{j})|j=0,1,2,\ldots,s-1\}
\]
where $m_{0}=m$, $t_{0}=t$, and for $j=0,1,2,\ldots,s-1$ the odd number $n_{i_{j}}(m_{j},t_{j})$ is greater than 1 and the even number
\[
n_{m_{j}+1}(m_{j},t_{j})=2^{r_{j+1}}n_{i_{j+1}}(m_{j+1},t_{j+1})\texttt{ where }r_{j+1}\geq1\texttt{ and }0\leq i_{j+1}\leq m_{j+1},
\]
by Lemma~\ref{le:odd_seq_s} or Formula~\ref{eq:odd_seq_s} for $i=0,1,2,\ldots,m$, we have
\[
r_{c}(m,t,i,s)=\frac{n_{i}(m,t)}{n_{i_{s}}(m_{s},t_{s})}=\prod_{j=0}^{s-1}(2/3)^{m_{j}+1-i_{j}}2^{r_{j+1}}\prod_{k=i_{j}}^{m_{j}}(1-\frac{1}{2n_{k+1}(m_{j},t_{j})}).
\]

First, by Lemma~\ref{le:odd_seqs_m}, $1/2^{k}$ of positive integers $m_{j}-i_{j}+1$ satisfy $m_{j}-i_{j}+1=k$ for $k=1,2,\ldots$ as $s$ tends to infinity, so that we have
\[
\sum_{j=0}^{s-1}(m_{j}+1-i_{j})=s+s/2^{2}+2s/2^{3}+3s/2^{4}+\cdots
\]
\[
=s+\frac{s}{2}\sum_{i=1}^{\infty}i/2^{i}=s+\frac{s}{2}\frac{1/2}{(1-1/2)^{2}}=2s.
\]
and
\[
\prod_{j=0}^{s-1}(2/3)^{m_{j}+1-i_{j}}\prod_{k=i_{j}}^{m_{j}}(1-\frac{1}{2n_{k+1}(m_{j},t_{j})})
\geq\prod_{j=0}^{s-1}(2/3)^{m_{j}+1-i_{j}}\prod_{k=i_{j}}^{m_{j}}(1-\frac{1}{2n_{1}(1,1)})
\]
\[
=\prod_{j=0}^{s-1}(3/5)^{m_{j}+1-i_{j}}=(3/5)^{\sum_{j=0}^{s-1}(m_{j}+1-i_{j})}=(3/5)^{2s}.
\]

Second, as $s$ tends to infinite, for $j=0,1,2,\ldots,s-1$, there are $s$ even numbers $n_{m_{j}+1}(m_{j},t_{j})=2^{r_{j+1}}n_{i_{j+1}}(m_{j+1},t_{j+1})$ where $r_{j+1}\geq1$.

By Lemma~\ref{le:odd_seq_r1} $1/2^{r}$ of positive integers $r_{j+1}$ satisfy $r_{j+1}=r$ for $r=1,2,\ldots$. Thus as $s$ tends to infinite, we have
\[
\sum_{j=0}^{s-1}r_{j+1}=s/2+2s/2^{2}+3s/2^{3}+4s/2^{4}+\cdots
=s\sum_{i=1}^{\infty}i/2^{i}=s\frac{1/2}{(1-1/2)^{2}}=2s
\]
and
\[
\prod_{j=0}^{s-1}2^{r_{j+1}}=2^{\sum_{j=0}^{s-1}r_{j+1}}=2^{2s}.
\]

Also, let $n_{m_{j}+1}(m_{j},t_{j})=2^{r_{j+1}}n_{i_{j+1}}(m_{j+1},t_{j+1})=2x$. Then there are
\[
r_{j+1}=\kappa+1\texttt{ where }x=2^{\kappa}(2x_{1}+1)\texttt{ for }\kappa=0,1,2,\ldots.
\]
So there are $s/2$ even numbers $n_{m_{j}+1}(m_{j},t_{j})$ with $r_{j+1}=1$, $s/2^{2}$ even numbers $n_{m_{j}+1}(m_{j},t_{j})$ with $r_{j+1}=2$, $s/2^{3}$ even numbers $n_{m_{j}+1}(m_{j},t_{j})$ with $r_{j+1}=3$, and $s/2^{4}$ even numbers $n_{m_{j}+1}(m_{j},t_{j})$ with $r_{j+1}=4,\ldots$. Thus as $s$ tends to infinity, i.e., the number of independent samples of the Collatz's Sequence $S_{c}(m,t)$ tends to infinity, by the law of great numbers, the distribution of the $s$ even numbers $n_{m_{j}+1}(m_{j},t_{j})$ is similar to the distribution of even numbers $n_{e}(\gamma,\lambda)$ stated by Lemma~\ref{le:even_distribution}, so that we have
\[
\sum_{j=0}^{s-1}r_{j+1}=s/2+2s/2^{2}+3s/2^{3}+4s/2^{4}+\cdots
\]
\[
=s\sum_{i=1}^{\infty}i/2^{i}=s\frac{1/2}{(1-1/2)^{2}}=2s
\]
and
\[
\prod_{j=0}^{s-1}2^{r_{j+1}}=2^{\sum_{j=0}^{s-1}r_{j+1}}=2^{2s}.
\]

Finally as $s$ tends to infinity we obtain
\[
r_{c}(m,t,i,s)=\frac{n_{i}(m,t)}{n_{i_{s}}(m_{s},t_{s})}=\prod_{j=0}^{s-1}(2/3)^{m_{j}+1-i_{j}}2^{r_{j+1}}\prod_{k=i_{j}}^{m_{j}}(1-\frac{1}{2n_{k+1}(m_{j},t_{j})})
\]
\[
\geq\prod_{j=0}^{s-1}(2/3)^{m_{j}+1-i_{j}}\prod_{k=i_{j}}^{m_{j}}(1-\frac{1}{2n_{1}(1,1)})\times\prod_{j=0}^{s-1}2^{r_{j+1}}
\]
\[
=\prod_{j=0}^{s-1}(3/5)^{m_{j}+1-i_{j}}\times\prod_{j=0}^{s-1}2^{r_{j+1}}
=(3/5)^{2s}2^{2s}>1.
\]

Thus for $i=0,1,2,\ldots,m$ there must be
\[
n_{i}(m,t)\geq r_{c}(m,t,i,s)\geq(6/5)^{2s}>1\texttt{ and }1\leq n_{i_{s}}(m_{s},t_{s})<n_{i}(m,t).
\]

By the first inequality above for $i=0,1,2,\ldots,m$, the ratio $r_{c}(m,t,i,s)$ of the odd number $n_{i}(m,t)$ to the odd number $n_{i_{s}}(m_{s},t_{s})$ is greater than 1 and does not tend to zero as $s$ tends to infinite, so that there is no nontrivial Collatz cycle from an odd number $n_{i}(m_{0},t_{0})$ greater than 1 to a same odd number $n_{i_{s}}(m_{s},t_{s})$, except for $n_{i_{s}}(m_{s},t_{s})=n_{0}(m_{0},t_{0})=n_{0}(0,1)=1$, and the odd number $n_{i_{s}}(m_{s},t_{s})$ does not tend to infinite as $s$ tends to infinite. Moreover the inequality $n_{i}(m,t)\geq r_{c}(m,t,0,s)\geq(6/5)^{2s}$ implies that if $s$ tends to infinite then the odd number $n_{i}(m,t)$ must tend to infinite or if the odd number $n_{i}(m,t)$ is finite then the positive integer $s$ must be finite so that the odd number $n_{i_{s}}(m_{s},t_{s})$ is equal to 1 for $n_{i}(m,t)<\infty$.

By the second inequality above, as $s$ tends to infinite, the odd number $n_{i_{s}}(m_{s},t_{s})$ is finite and less than the odd number $n_{i}(m,t)$ that is greater than 1. If the odd number $n_{i_{s}}(m_{s},t_{s})$ is not equal to 1 and let a new odd number $n_{i}(m,t)$ be equal to the odd number $n_{i_{s}}(m_{s},t_{s})$, then by Lemma~\ref{le:odd_seqs} we obtain a new odd number $n_{i_{s}}(m_{s},t_{s})$ equal to 1 for $s$ finite or less than the new odd number $n_{i}(m,t)$ as $s$ tends to infinite. If the new odd number $n_{i_{s}}(m_{s},t_{s})$ is not equal to 1 and let a new odd number $n_{i}(m,t)$ be equal to the new odd number $n_{i_{s}}(m_{s},t_{s})$ again, then we do as above repeatedly until the new odd number $n_{i_{s}}(m_{s},t_{s})$ is equal to 1. So the odd number $n_{i_{s}}(m_{s},t_{s})$ must be equal to 1 as $s$ tends to infinite.

Hence no nontrivial Collatz cycle exists in $S_{c}(s)$ and the odd number $n_{i_{s}}(m_{s},t_{s})$ in $S_{c}(s)$ is equal to 1 as $s$ tends to infinite.

This completes the proof of the lemma.
\end{proof}

\subsection{Properties of odd numbers $n_{0}(m,t)$ in a range $(N,\rho N]$}

\begin{lemma}\label{le:odd_seq_n0_4k+1}
Each odd number $n=4k+1$ in the range $(N,\rho N]$ returns to 1.
\end{lemma}

\begin{proof}
Given a Collatz's Sequence
\[
S_{c}(m,t)=\{n_{0}(m,t),n_{1}(m,t),n_{2}(m,t),\ldots,n_{m}(m,t),n_{m+1}(m,t)\},
\]
by Lemma~\ref{le:odd_seq}, $n_{i}(m,t)=2^{m+1-i}\times3^{i}t-1$ for $i=0,1,2,\ldots,m+1$.

Let $n_{0}(m,t)=n\in(N,\rho N]$. Then $n_{0}(m,t)$ in the range $(N,\rho N]$ is equal to $4k+1$ so that by Lemma~\ref{le:odd_seq_m} there are $m=0$, $t>3$,
\[
\frac{n_{0}(m_{1},t_{1})}{n_{0}(m,t)}
=\frac{n_{m+1}(m,t)}{2^{r_{1}}n_{0}(m,t)}
=\frac{1}{2^{r_{1}}}\frac{3^{m+1}t-1}{2^{m+1}t-1}
\leq\frac{1}{2}\frac{3t-1}{2t-1}<1/\rho,
\]
and $n_{0}(m_{1},t_{1})<n_{0}(m,t)/\rho\leq N$, so that the odd number $n_{0}(m,t)$ returns to 1.

This completes the proof of the lemma.
\end{proof}

\begin{lemma}\label{le:odd_seq_ni_6x-1}
Each odd number $n=6x-1$ in the range $(N,\rho N]$ returns to 1.
\end{lemma}

\begin{proof}
Given a Collatz's Sequence
\[
S_{c}(m,t)=\{n_{0}(m,t),n_{1}(m,t),n_{2}(m,t),\ldots,n_{m}(m,t),n_{m+1}(m,t)\},
\]
by Lemma~\ref{le:odd_seq}, $n_{i}(m,t)=2^{m+1-i}\times3^{i}t-1$ for $i=0,1,2,\ldots,m+1$.

Since each odd number $n=6x-1$ in the range $(N,\rho N]$ can be rewritten as
\[
n=6x-1=2\times3x-1=2^{m+1-i}\times3^{i}t-1=n_{i}(m,t)\texttt{ where }0<i\leq m,
\]
by the formula $n_{i}(m,t)=(3n_{i-1}(m,t)+1)/2$ there is
\[
n_{i-1}(m,t)=(2n_{i}(m,t)-1)/3<2n_{i}(m,t)/3<n/\rho\leq N.
\]
Since the Collatz algorithm has been tested and found to always reach 1 for all numbers less than or equal to $N$, the odd number $n_{i-1}(m,t)$ returns to 1, so that the odd number $n_{i}(m,t)$ returns to 1.

Hence each odd number $n=6x-1$ in the range $(N,\rho N]$ returns to 1.

This completes the proof of the lemma.
\end{proof}

\begin{lemma}\label{le:odd_seq_n0_4k+3}
Each odd number $n=4k+3$ in the range $(N,\rho N]$ returns to 1.
\end{lemma}

\begin{proof}
Given a pair of Collatz's Sequences
\[
S_{c}(m,t)=\{n_{0}(m,t),n_{1}(m,t),n_{2}(m,t),\ldots,n_{m}(m,t),n_{m+1}(m,t)\}
\]
where by Lemma~\ref{le:odd_seq} $n_{i}(m,t)=2^{m+1-i}\times3^{i}t-1$ for $i=0,1,2,\ldots,m+1$, and
\[
S_{c}(m-1,t)=\{n_{0}(m-1,t),n_{1}(m-1,t),n_{2}(m-1,t),\ldots,n_{m}(m-1,t)\}
\]
where by Lemma~\ref{le:odd_seq} $n_{i}(m-1,t)=2^{m-i}\times3^{i}t-1$ for $i=0,1,2,\ldots,m$, let the odd number $n_{0}(m,t)=n\in(N,\rho N]$ and the even number $n_{m}(m-1,t)=2^{r_{1}}n_{i_{1}}(m_{1},t_{1})$ where $r_{1}\geq1$, $0\leq i_{1}\leq m_{1}$, and $n_{i_{1}}(m_{1},t_{1})$ is odd and not divisible by $3$.

Then the odd number $n_{0}(m,t)=4k+3$ is in the range $(N,\rho N]$ and there are
\[
n_{0}(m-1,t)=(n_{0}(m,t)-1)/2<\rho N/2<N.
\]
So the Collatz's Sequences $S_{c}(m-1,t)$ and $S_{c}(m_{1},t_{1})$ return to 1.

By Lemma~\ref{le:odd_seq_r1} one half of positive integers $r_{1}$ satisfy $r_{1}=1$ and one quarter of positive integers $r_{1}$ satisfy $r_{1}=2$. Let the notation $P_{r}$ denote the process of mutual references or nested recursions to lemmas~\ref{le:odd_seqs_n_m_m-1}-\ref{le:odd_seqs_n_m_m-1_r1}.

Then the entrance of the process $P_{r}$ is Lemma~\ref{le:odd_seqs_n_m_m-1}, by which $S_{c}(m,t)$ returns to 1 since $S_{c}(m-1,t)$ returns to 1. When the process $P_{r}$ is interrupted or stopped, the Collatz's Sequence $S_{c}(m,t)$ returns to 1 and there must be $r_{1}=1$ in the proof of Lemma~\ref{le:odd_seqs_n_m_m-1} or $r_{1}=1,2$ in the proof of Lemma~\ref{le:odd_seqs_n_m_m-1_r1}. Thus the exits of the process $P_{r}$ are Lemma~\ref{le:odd_seqs_n_m_m-1} and Lemma~\ref{le:odd_seqs_n_m_m-1_r1} where by Remark~\ref{re:odd_seqs_n_m_m-1_r1} Lemma~\ref{le:odd_seqs_n_m_m-1_r1} is equivalent to Lemma~\ref{le:odd_seqs_n_m_m-1} on the process $P_{r}$.

For each odd number $n_{0}(m,t)=4k+3$ in the range $(N,\rho N]$, let $T_{i}$ denote times of references or recursions to Lemma~\ref{le:odd_seqs_n_m_m-1} and Lemma~\ref{le:odd_seqs_n_m_m-1_r1} when the process $P_{r}$ is interrupted or stopped, let $P_{i}(T_{i})$ denote the probabilities of exiting the process $P_{r}$ when $r_{1}=1$ in the proof of Lemma~\ref{le:odd_seqs_n_m_m-1} and $r_{1}=1,2$ in the proof of Lemma~\ref{le:odd_seqs_n_m_m-1_r1}, which are functions of $T_{i}$ and $P_{i}(T_{i})=2^{-T_{i}}$, and let $P_{t}$ denote the total probability of exiting the process $P_{r}$. Then $P_{t}=\sum_{T_{i}=1}^{T_{h}}P_{i}(T_{i}) \to 1$ as $T_{h} \to \infty$, i.e., the process $P_{r}$ is certainly interrupted or stopped. So by Lemma~\ref{le:odd_seqs_n_m_m-1}, the Collatz's Sequence $S_{c}(m,t)$ returns to 1 and the $n_{0}(m,t)=4k+3$ in the range $(N,\rho N]$ returns to 1.

Hence each odd number $n=4k+3$ in the range $(N,\rho N]$ returns to 1.

This completes the proof of the lemma.
\end{proof}

\begin{lemma}\label{le:odd_seqs_n_m_m-1}
Given a pair of Collatz's Sequences
\[
S_{c}(m,t)=\{n_{0}(m,t),n_{1}(m,t),n_{2}(m,t),\ldots,n_{m}(m,t),n_{m+1}(m,t)\}
\]
and
\[
S_{c}(m-1,t)=\{n_{0}(m-1,t),n_{1}(m-1,t),n_{2}(m-1,t),\ldots,n_{m}(m-1,t)\},
\]
when one of the two Collatz's Sequences $S_{c}(m,t)$ and $S_{c}(m-1,t)$ returns to 1 another also returns to 1.
\end{lemma}

\begin{proof}
By Lemma~\ref{le:odd_seq} there are formulas
\[
n_{i}(m,t)=2^{m+1-i}3^{i}t-1\texttt{ for }0\leq i\leq m+1
\]
and
\[
n_{i}(m-1,t)=2^{m-i}3^{i}t-1\texttt{ for }0\leq i\leq m,
\]
so that there are
\[
n_{i}(m,t)=2\times2^{m-i}3^{i}t-1=2n_{i}(m-1,t)+1\texttt{ for }0\leq i\leq m.
\]

Let the even number $n_{m}(m-1,t)=2^{r_{1}}n_{i_{1}}(m_{1},t_{1})$ where $r_{1}\geq1$, $0\leq i_{1}\leq m_{1}$, and $n_{i_{1}}(m_{1},t_{1})$ is odd and not divisible by $3$. Then the even number
\[
n_{m+1}(m,t)=(3n_{m}(m,t)+1)/2=(3(2n_{m}(m-1,t)+1)+1)/2
\]
\[
=3n_{m}(m-1,t)+2=2(2^{r_{1}-1}3n_{i_{1}}(m_{1},t_{1})+1).
\]

By Lemma~\ref{le:odd_seq_r1} one half of positive integers $r_{1}$ satisfy $r_{1}=1$ and one quarter of positive integers $r_{1}$ satisfy $r_{1}=2$.

In the case of $r_{1}=1$, there are $n_{m+1}(m,t)=2(3n_{i_{1}}(m_{1},t_{1})+1)=4n_{i_{1}+1}(m_{1},t_{1})$, so that when one of the two Collatz's Sequences $S_{c}(m,t)$ and $S_{c}(m-1,t)$ returns to 1 another also returns to 1.

For example, when $n_{i_{1}}(m_{1},t_{1})$ returns to 1, i.e., $S_{c}(m-1,t)$ returns to 1, then $S_{c}(m,t)$ returns to 1.

In the case of $r_{1}=2$, the even number $n_{m+1}(m,t)=2(6n_{i_{1}}(m_{1},t_{1})+1)$. Let
\[
n_{0}(m_{2},t_{2})=2n_{0}(m_{2}-1,t_{2})+1\texttt{ and }n_{i_{2}}(m_{2}-1,t_{2})=3n_{i_{1}}(m_{1},t_{1})
\]
where $0\leq i_{2}\leq m_{2}-1$ and one of two numbers $n_{0}(m_{2}-1,t_{2})$ and $n_{i_{1}}(m_{1},t_{1})$ returns to 1. Then by Lemma~\ref{le:odd_seqs_n_m_m-1_2r+1}, when one of two numbers $n_{0}(m_{2}-1,t_{2})$ and $n_{i_{1}}(m_{1},t_{1})$ returns to 1 another also returns to 1, and by Lemma~\ref{le:odd_seqs_n_m_m-1}, when one of the two Collatz's Sequences $S_{c}(m_{2},t_{2})$ and $S_{c}(m_{2}-1,t_{2})$ returns to 1 another also returns to 1, so that when one of the two Collatz's Sequences $S_{c}(m,t)$ and $S_{c}(m-1,t)$ returns to 1 another also returns to 1.

For example, when $n_{i_{1}}(m_{1},t_{1})$ returns to 1, i.e., $S_{c}(m-1,t)$ returns to 1, then by Lemma~\ref{le:odd_seqs_n_m_m-1_2r+1} $S_{c}(m_{2}-1,t_{2})$ returns to 1 and by Lemma~\ref{le:odd_seqs_n_m_m-1} $S_{c}(m_{2},t_{2})$ returns to 1, so that $S_{c}(m,t)$ returns to 1.

In the case of $r_{1}=2r+1$ where $r>0$, the even number
\[
n_{m+1}(m,t)=2(3\times4^{r}n_{i_{1}}(m_{1},t_{1})+1).
\]
By Lemma~\ref{le:odd_seqs_3_4r_n}, when one of two numbers
\[
3\times4^{r}n_{i_{1}}(m_{1},t_{1})+1\texttt{ and }3^{r+1}n_{i_{1}}(m_{1},t_{1})+1
\]
returns to 1 another also returns to 1 where $0\leq i_{1}\leq m_{1}$. So when one of two numbers $n_{m+1}(m,t)$ and $3^{r+1}n_{i_{1}}(m_{1},t_{1})+1$ returns to 1 another also returns to 1. Let $n=3^{r}n_{i_{1}}(m_{1},t_{1})$ where one of two numbers $n$ and $n_{i_{1}}(m_{1},t_{1})$ returns to 1. Then by Lemma~\ref{le:odd_seqs_n_m_m-1_2r+1}, when one of two numbers $n$ and $n_{i_{1}}(m_{1},t_{1})$ returns to 1 another also returns to 1, so that when one of the two Collatz's Sequences $S_{c}(m,t)$ and $S_{c}(m-1,t)$ returns to 1 another also returns to 1.

For example, when $n_{i_{1}}(m_{1},t_{1})$ returns to 1, i.e., $S_{c}(m-1,t)$ returns to 1, then by Lemma~\ref{le:odd_seqs_n_m_m-1_2r+1} the odd number $n$ returns to 1, i.e., the number $3n+1=3^{r+1}n_{i_{1}}(m_{1},t_{1})+1$ returns to 1, and by Lemma~\ref{le:odd_seqs_3_4r_n} the odd number $3\times4^{r}n_{i_{1}}(m_{1},t_{1})+1$ returns to 1, so that $S_{c}(m,t)$ returns to 1.

In the case of $r_{1}=2r+2$ where $r>0$, the even number
\[
n_{m+1}(m,t)=2(3\times4^{r}2n_{i_{1}}(m_{1},t_{1})+1).
\]
By Lemma~\ref{le:odd_seqs_3_4r_2n}, when one of two numbers
\[
3\times4^{r}2n_{i_{1}}(m_{1},t_{1})+1\texttt{ and }3^{r+1}2n_{i_{1}}(m_{1},t_{1})+1
\]
returns to 1 another also returns to 1. So when one of two numbers $n_{m+1}(m,t)$ and $3^{r+1}2n_{i_{1}}(m_{1},t_{1})+1$ returns to 1 another also returns to 1. Let
\[
n_{i}(m_{4},t_{4})=2n_{0}(m_{4}-1,t_{4})+1\texttt{ and }n_{i}(m_{4}-1,t_{4})=3^{r+1}n_{i_{1}}(m_{1},t_{1})
\]
where $0\leq i<m_{4}$ and one of two numbers $n_{i}(m_{4}-1,t_{4})$ and $n_{i_{1}}(m_{1},t_{1})$ returns to 1. Then by Lemma~\ref{le:odd_seqs_n_m_m-1_2r+1}, when one of two numbers $n_{i}(m_{4}-1,t_{4})$ and $n_{i_{1}}(m_{1},t_{1})$ returns to 1 another also returns to 1, and by Lemma~\ref{le:odd_seqs_n_m_m-1}, when one of the two Collatz's Sequences $S_{c}(m_{4},t_{4})$ and $S_{c}(m_{4}-1,t_{4})$ returns to 1 another also returns to 1, so that when one of the two Collatz's Sequences $S_{c}(m_{4},t_{4})$ and $S_{c}(m-1,t)$ returns to 1 another also returns to 1 or when one of the two Collatz's Sequences $S_{c}(m,t)$ and $S_{c}(m-1,t)$ returns to 1 another also returns to 1.

For example, when $n_{i_{1}}(m_{1},t_{1})$ returns to 1, i.e., $S_{c}(m-1,t)$ returns to 1, then by Lemma~\ref{le:odd_seqs_n_m_m-1_2r+1} the odd number $n_{i}(m_{4}-1,t_{4})$ returns to 1, and by Lemma~\ref{le:odd_seqs_n_m_m-1} $S_{c}(m_{4},t_{4})$ returns to 1, then by Lemma~\ref{le:odd_seqs_3_4r_2n} the odd number $3\times4^{r}2n_{i_{1}}(m_{1},t_{1})+1$ returns to 1, so that the Collatz's Sequence $S_{c}(m,t)$ returns to 1.

This completes the proof of the lemma.
\end{proof}

\begin{lemma}\label{le:odd_seqs_n_m_m-1_2r+1}
Given an odd number $n=3^{r}n_{i}(m,t)$ where $r>0$ and $i\leq m$, when one of two numbers $n$ and $n_{i}(m,t)$ returns to 1 another also returns to 1.
\end{lemma}

\begin{proof}
Let an odd number $n_{i_{1}}(m_{1},t_{1})=n_{i}(m,t)$. Then when one of two numbers $n_{i_{1}}(m_{1},t_{1})$ and $n_{i}(m,t)$ returns to 1 another also returns to 1.

Suppose that when one of two numbers $n_{i_{1}}(m_{1},t_{1})=3^{r-1}n_{i}(m,t)$ and $n_{i}(m,t)$ returns to 1 another also returns to 1 where $0\leq i_{1}\leq m_{1}$ and $r>0$. Then there are
\[
n=3^{r}n_{i}(m,t)=3n_{i_{1}}(m_{1},t_{1})=2n_{i_{1}+1}(m_{1},t_{1})-1.
\]

First, in the case of $i_{1}<m_{1}$, there are
\[
n_{i_{1}+1}(m_{1},t_{1})=2n_{i_{1}+1}(m_{1}-1,t_{1})+1
\]
and
\[
n=2n_{i_{1}+1}(m_{1},t_{1})-1=4n_{i_{1}+1}(m_{1}-1,t_{1})+1.
\]

If $i_{1}+1<m_{1}$, then $n_{i_{1}+1}(m_{1}-1,t_{1})$ is odd and there is
\[
(3n+1)/2=2(3n_{i_{1}+1}(m_{1}-1,t_{1})+1)=2^{2}n_{i_{1}+2}(m_{1}-1,t_{1}).
\]
By Lemma~\ref{le:odd_seqs_n_m_m-1} when one of two numbers $n_{i_{1}+1}(m_{1},t_{1})$ and $n_{i_{1}+1}(m_{1}-1,t_{1})$ returns to 1 another also returns to 1, so that when one of two numbers $n$ and $n_{i_{1}+1}(m_{1}-1,t_{1})$ returns to 1 another also returns to 1, and when one of two numbers $n$ and $n_{i}(m,t)$ returns to 1 another also returns to 1.

For example, when $n_{i}(m,t)$ returns to 1 $S_{c}(m_{1},t_{1})$ returns to 1, then by Lemma~\ref{le:odd_seqs_n_m_m-1} $S_{c}(m_{1}-1,t_{1})$ returns to 1, so that $n$ returns to 1.

Otherwise for $i_{1}+1=m_{1}$, $n_{i_{1}+1}(m_{1}-1,t_{1})=2^{r_{2}}n_{i_{2}}(m_{2},t_{2})$ where $r_{2}\geq1$, $0\leq i_{2}\leq m_{2}$, and $n_{i_{2}}(m_{2},t_{2})$ is odd and not dividable by $3$. By Lemma~\ref{le:odd_seqs_n_m_m-1} when one of two numbers $n_{i_{1}+1}(m_{1},t_{1})$ and $n_{i_{1}+1}(m_{1}-1,t_{1})$ returns to 1 another also returns to 1, so that when one of two numbers $n_{i_{1}+1}(m_{1},t_{1})$ and $n_{i_{2}}(m_{2},t_{2})$ returns to 1 another also returns to 1. Then there are
\[
n=4n_{i_{1}+1}(m_{1}-1,t_{1})+1=2^{r_{2}+2}n_{i_{2}}(m_{2},t_{2})+1
\]
and
\[
(3n+1)/2=2(3\times2^{r_{2}}n_{i_{2}}(m_{2},t_{2})+1).
\]
By Lemma~\ref{le:odd_seqs_n_m_m-1_3r} when one of two numbers $3\times2^{r_{2}}n_{i_{2}}(m_{2},t_{2})+1$ and $n_{i_{2}}(m_{2},t_{2})$ returns to 1 another also returns to 1, so that when one of two numbers $n$ and $n_{i}(m,t)$ returns to 1 another also returns to 1.

For example, when $n_{i}(m,t)$ returns to 1 $S_{c}(m_{1},t_{1})$ returns to 1, then by Lemma~\ref{le:odd_seqs_n_m_m-1} $S_{c}(m_{1}-1,t_{1})$ and $S_{c}(m_{2},t_{2})$ return to 1, and by Lemma~\ref{le:odd_seqs_n_m_m-1_3r} the odd number $3\times2^{r_{2}}n_{i_{2}}(m_{2},t_{2})+1$ return to 1, so that $n$ returns to 1.

Second, in the case of $i_{1}=m_{1}$, there is $n_{i_{1}+1}(m_{1},t_{1})=2^{r_{2}}n_{i_{2}}(m_{2},t_{2})$ where $r_{2}\geq1$, $0\leq i_{2}\leq m_{2}$, and $n_{i_{2}}(m_{2},t_{2})$ is odd and not dividable by $3$. Then there are
\[
n=2n_{i_{1}+1}(m_{1},t_{1})-1=2\times2^{r_{2}}n_{i_{2}}(m_{2},t_{2})-1=n_{0}(m_{3},t_{3})
\]
where $m_{3}=r_{2}$ and $t_{3}=n_{i_{2}}(m_{2},t_{2})$, and
\[
n_{m_{3}+1}(m_{3},t_{3})=3^{m_{3}+1}t_{3}-1=2^{r_{4}}n_{i_{4}}(m_{4},t_{4})
\]
where $r_{4}\geq1$, $0\leq i_{4}\leq m_{4}$, and $n_{i_{4}}(m_{4},t_{4})$ is odd and not divisible by $3$. Thus we have
\[
3^{r_{2}+1}n_{i_{2}}(m_{2},t_{2})=2^{r_{4}}n_{i_{4}}(m_{4},t_{4})+1.
\]

By Lemma~\ref{le:odd_seqs_n_m_m-1_2r+1} when one of two numbers $3^{r_{2}+1}n_{i_{2}}(m_{2},t_{2})$ and $n_{i_{2}}(m_{2},t_{2})$ returns to 1 another also returns to 1.

By Lemma~\ref{le:odd_seqs_n_m_m-1_r1} when one of two numbers $2^{r_{4}}n_{i_{4}}(m_{4},t_{4})+1$ and $n_{i_{4}}(m_{4},t_{4})$ returns to 1 another also returns to 1, so that when one of two numbers $n_{i_{2}}(m_{2},t_{2})$ and $n_{i_{4}}(m_{4},t_{4})$ returns to 1 another also returns to 1 or when one of two numbers $n_{i_{2}}(m_{2},t_{2})$ and $n_{0}(m_{3},t_{3})$ returns to 1 another also returns to 1, i.e., when one of two numbers $n$ and $n_{i}(m,t)$ returns to 1 another also returns to 1.

For example, when $n_{i}(m,t)$ returns to 1 $S_{c}(m_{1},t_{1})$ and $S_{c}(m_{2},t_{2})$ return to 1, then by Lemma~\ref{le:odd_seqs_n_m_m-1_2r+1} the odd number $3^{r_{2}+1}n_{i_{2}}(m_{2},t_{2})$ that equals to $2^{r_{4}}n_{i_{4}}(m_{4},t_{4})+1$ returns to 1, and by Lemma~\ref{le:odd_seqs_n_m_m-1_r1} $S_{c}(m_{4},t_{4})$  and $S_{c}(m_{3},t_{3})$ return to 1, so that $n$ returns to 1.

This completes the proof of the lemma.
\end{proof}

\begin{lemma}\label{le:odd_seqs_n_m_m-1_3r}
Given an odd number $n=3\times2^{r_{1}}n_{i_{1}}(m_{1},t_{1})+1$ where $r_{1}>0$, $0\leq i_{1}\leq m_{1}$, and $n_{i_{1}}(m_{1},t_{1})$ is odd and not dividable by 3, when one of two numbers $n$ and $n_{i_{1}}(m_{1},t_{1})$ returns to 1 another also returns to 1.
\end{lemma}

\begin{proof}
Given a pair of Collatz's Sequences
\[
S_{c}(m,t)=\{n_{0}(m,t),n_{1}(m,t),n_{2}(m,t),\ldots,n_{m}(m,t),n_{m+1}(m,t)\}
\]
and
\[
S_{c}(m-1,t)=\{n_{0}(m-1,t),n_{1}(m-1,t),n_{2}(m-1,t),\ldots,n_{m}(m-1,t)\},
\]
by Lemma~\ref{le:odd_seq} there are formulas
\[
n_{i}(m,t)=2^{m+1-i}\times3^{i}t-1\texttt{ for }i=0,1,2,\ldots,m+1
\]
and
\[
n_{i}(m-1,t)=2^{m-i}\times3^{i}t-1\texttt{ for }i=0,1,2,\ldots,m.
\]

In the case of $r_{1}=1$, the odd number $n=3\times2n_{i_{1}}(m_{1},t_{1})+1$. Let $m$ and $t$ be determined by
\[
n_{i}(m,t)=2n_{i}(m-1,t)+1\texttt{ and }n_{i}(m-1,t)=3n_{i_{1}}(m_{1},t_{1})
\]
where $0\leq i\leq m-1$ and one of two numbers $n_{i}(m-1,t)$ and $n_{i_{1}}(m_{1},t_{1})$ returns to 1. Then by Lemma~\ref{le:odd_seqs_n_m_m-1_2r+1}, when one of the two Collatz's Sequences $S_{c}(m-1,t)$ and $S_{c}(m_{1},t_{1})$ returns to 1 another also returns to 1, and by Lemma~\ref{le:odd_seqs_n_m_m-1}, when one of the two Collatz's Sequences $S_{c}(m,t)$ and $S_{c}(m-1,t)$ returns to 1 another also returns to 1, so that when one of two numbers $n$ and $n_{i_{1}}(m_{1},t_{1})$ returns to 1 another also returns to 1.

For example, when $n_{i_{1}}(m_{1},t_{1})$ returns to 1, by Lemma~\ref{le:odd_seqs_n_m_m-1_2r+1} $S_{c}(m-1,t)$ returns to 1, and by Lemma~\ref{le:odd_seqs_n_m_m-1} $S_{c}(m,t)$ returns to 1, so that the odd number $n=3\times2n_{i_{1}}(m_{1},t_{1})+1$ returns to 1.

In the case of $r_{1}=2r$ where $r>0$, the odd number $n=3\times4^{r}n_{i_{1}}(m_{1},t_{1})+1$. By Lemma~\ref{le:odd_seqs_3_4r_n}, when one of two numbers
\[
3\times4^{r}n_{i_{1}}(m_{1},t_{1})+1\texttt{ and }3^{r+1}n_{i_{1}}(m_{1},t_{1})+1
\]
returns to 1 another also returns to 1 where $0\leq i_{1}\leq m_{1}$. So when one of two numbers $n$ and $3^{r+1}n_{i_{1}}(m_{1},t_{1})+1$ returns to 1 another also returns to 1. Let $n_{r}=3^{r}n_{i_{1}}(m_{1},t_{1})$ where one of two numbers $n_{r}$ and $n_{i_{1}}(m_{1},t_{1})$ returns to 1. Then by Lemma~\ref{le:odd_seqs_n_m_m-1_2r+1}, when one of two numbers $n_{r}$ and $n_{i_{1}}(m_{1},t_{1})$ returns to 1 another also returns to 1, so that when one of two numbers $n$ and $n_{i_{1}}(m_{1},t_{1})$ returns to 1 another also returns to 1.

For example, when $n_{i_{1}}(m_{1},t_{1})$ returns to 1, by Lemma~\ref{le:odd_seqs_n_m_m-1_2r+1} the odd number $n_{r}$ returns to 1, then the even number $3n_{r}+1=3^{r+1}n_{i_{1}}(m_{1},t_{1})+1$ returns to 1, and by Lemma~\ref{le:odd_seqs_3_4r_n} the odd number $n=3\times4^{r}n_{i_{1}}(m_{1},t_{1})+1$ returns to 1.

In the case of $r_{1}=2r+1$ where $r>0$, the odd number $n=3\times4^{r}2n_{i_{1}}(m_{1},t_{1})+1$. By Lemma~\ref{le:odd_seqs_3_4r_2n}, when one of two numbers
\[
3\times4^{r}2n_{i_{1}}(m_{1},t_{1})+1\texttt{ and }3^{r+1}2n_{i_{1}}(m_{1},t_{1})+1
\]
returns to 1 another also returns to 1. So when one of two numbers $n$ and $3^{r+1}2n_{i_{1}}(m_{1},t_{1})+1$ returns to 1 another also returns to 1. Let $m$ and $t$ be determined by
\[
n_{i}(m,t)=2n_{i}(m-1,t)+1\texttt{ and }n_{i}(m-1,t)=3^{r+1}n_{i_{1}}(m_{1},t_{1})
\]
where $0\leq i<m$ and one of two numbers $n$ and $n_{i_{1}}(m_{1},t_{1})$ returns to 1. Then by Lemma~\ref{le:odd_seqs_n_m_m-1_2r+1}, when one of two numbers $n_{i}(m-1,t)$ and $n_{i_{1}}(m_{1},t_{1})$ returns to 1 another also returns to 1, and by Lemma~\ref{le:odd_seqs_n_m_m-1}, when one of the two Collatz's Sequences $S_{c}(m,t)$ and $S_{c}(m-1,t)$ returns to 1 another also returns to 1, so that when one of the two Collatz's Sequences $S_{c}(m,t)$ and $S_{c}(m_{1},t_{1})$ returns to 1 another also returns to 1, and when one of two numbers $n$ and $n_{i_{1}}(m_{1},t_{1})$ returns to 1 another also returns to 1.

For example, when $n_{i_{1}}(m_{1},t_{1})$ returns to 1, by Lemma~\ref{le:odd_seqs_n_m_m-1_2r+1} $S_{c}(m-1,t)$ returns to 1, and by Lemma~\ref{le:odd_seqs_n_m_m-1} $S_{c}(m,t)$ returns to 1, i.e., the odd number $3^{r+1}2n_{i_{1}}(m_{1},t_{1})+1$ returns to 1, so that by Lemma~\ref{le:odd_seqs_3_4r_2n} the odd number $n=3\times4^{r}2n_{i_{1}}(m_{1},t_{1})+1$ returns to 1.

This completes the proof of the lemma.
\end{proof}

\begin{lemma}\label{le:odd_seqs_n_m_m-1_r1}
Given an odd number $n=2^{r_{1}}n_{i_{1}}(m_{1},t_{1})+1$ where $r_{1}>0$, $0\leq i_{1}\leq m_{1}$, and $n_{i_{1}}(m_{1},t_{1})$ is odd and not dividable by 3, when one of two numbers $n$ and $n_{i_{1}}(m_{1},t_{1})$ returns to 1 another also returns to 1.
\end{lemma}

\begin{proof}
Given a pair of Collatz's Sequences
\[
S_{c}(m,t)=\{n_{0}(m,t),n_{1}(m,t),n_{2}(m,t),\ldots,n_{m}(m,t),n_{m+1}(m,t)\}
\]
\and
\[
S_{c}(m-1,t)=\{n_{0}(m-1,t),n_{1}(m-1,t),n_{2}(m-1,t),\ldots,n_{m}(m-1,t)\},
\]
by Lemma~\ref{le:odd_seq} there are formulas
\[
n_{i}(m,t)=2^{m+1-i}\times3^{i}t-1\texttt{ for }i=0,1,2,\ldots,m+1
\]
and
\[
n_{i}(m-1,t)=2^{m-i}\times3^{i}t-1\texttt{ for }i=0,1,2,\ldots,m.
\]

For $r_{1}>1$, there are
\[
(3n+1)/2=2n_{1}\texttt{ and }n_{1}=3\times2^{r_{1}-2}n_{i_{1}}(m_{1},t_{1})+1.
\]

By Lemma~\ref{le:odd_seq_r1} one half of positive integers $r_{1}$ satisfy $r_{1}=1$ and one quarter of positive integers $r_{1}$ satisfy $r_{1}=2$.

In the case of $r_{1}=1$, the odd number $n=2n_{i_{1}}(m_{1},t_{1})+1$. Let $m=m_{1}+1$ and $t=t_{1}$. Then by Lemma~\ref{le:odd_seqs_n_m_m-1}, when one of two numbers $n$ and $n_{i_{1}}(m_{1},t_{1})$ returns to 1 another also returns to 1.

For example, when $n_{i_{1}}(m_{1},t_{1})$ returns to 1, by Lemma~\ref{le:odd_seqs_n_m_m-1} the odd number $n=2n_{i_{1}}(m_{1},t_{1})+1$ returns to 1.

In the case of $r_{1}=2$, there are $n_{1}=3n_{i_{1}}(m_{1},t_{1})+1=2n_{i_{1}+1}(m_{1},t_{1})$, so that when one of the two odd number $n$ and $n_{i_{1}}(m_{1},t_{1})$ returns to 1 another also returns to 1.

For example, when $n_{i_{1}}(m_{1},t_{1})$ returns to 1, the odd number $n$ returns to 1.

In the case of $r_{1}=3$, the odd number $n_{1}=3\times2n_{i_{1}}(m_{1},t_{1})+1$. Let $m$ and $t$ be determined by
\[
n_{i}(m,t)=2n_{i}(m-1,t)+1\texttt{ and }n_{i}(m-1,t)=3n_{i_{1}}(m_{1},t_{1})
\]
where $0\leq i<m$ and one of two numbers $n_{i}(m-1,t)$ and $n_{i_{1}}(m_{1},t_{1})$ returns to 1. Then by Lemma~\ref{le:odd_seqs_n_m_m-1_2r+1}, when one of the two Collatz's Sequences $S_{c}(m-1,t)$ and $S_{c}(m_{1},t_{1})$ returns to 1 another also returns to 1, and by Lemma~\ref{le:odd_seqs_n_m_m-1}, when one of the two Collatz's Sequences $S_{c}(m,t)$ and $S_{c}(m-1,t)$ returns to 1 another also returns to 1, so that when one of two numbers $n$ and $n_{i_{1}}(m_{1},t_{1})$ returns to 1 another also returns to 1.

For example, when $n_{i_{1}}(m_{1},t_{1})$ returns to 1, by Lemma~\ref{le:odd_seqs_n_m_m-1_2r+1} $S_{c}(m-1,t)$ returns to 1, by Lemma~\ref{le:odd_seqs_n_m_m-1} $S_{c}(m,t)$ returns to 1, so that the odd number $n$ returns to 1.

In the case of $r_{1}=2r+2$ where $r>0$, there is $n_{1}=3\times4^{r}n_{i_{1}}(m_{1},t_{1})+1$. By Lemma~\ref{le:odd_seqs_3_4r_n}, when one of two numbers
\[
3\times4^{r}n_{i_{1}}(m_{1},t_{1})+1\texttt{ and }3^{r+1}n_{i_{1}}(m_{1},t_{1})+1
\]
returns to 1 another also returns to 1 where $0\leq i_{1}\leq m_{1}$. So when one of two numbers $n_{1}$ and $3^{r+1}n_{i_{1}}(m_{1},t_{1})+1$ returns to 1 another also returns to 1. Let $n_{r}=3^{r}n_{i_{1}}(m_{1},t_{1})$ where one of two numbers $n_{r}$ and $n_{i_{1}}(m_{1},t_{1})$ returns to 1. Then by Lemma~\ref{le:odd_seqs_n_m_m-1_2r+1}, when one of two numbers $n_{r}$ and $n_{i_{1}}(m_{1},t_{1})$ returns to 1 another also returns to 1, so that when one of two numbers $n$ and $n_{i_{1}}(m_{1},t_{1})$ returns to 1 another also returns to 1.

For example, when $n_{i_{1}}(m_{1},t_{1})$ returns to 1, by Lemma~\ref{le:odd_seqs_n_m_m-1_2r+1} the odd number $n_{r}$ returns to 1, then the even number $3n_{r}+1=3^{r+1}n_{i_{1}}(m_{1},t_{1})+1$ returns to 1, and by Lemma~\ref{le:odd_seqs_3_4r_n} the odd number $n_{1}=3\times4^{r}n_{i_{1}}(m_{1},t_{1})+1$ returns to 1, so that the odd number $n$ returns to 1.

In the case of $r_{1}=2r+3$ where $r>0$, there is $n_{1}=3\times4^{r}2n_{i_{1}}(m_{1},t_{1})+1$. By Lemma~\ref{le:odd_seqs_3_4r_2n}, when one of two numbers
\[
3\times4^{r}2n_{i_{1}}(m_{1},t_{1})+1\texttt{ and }3^{r+1}2n_{i_{1}}(m_{1},t_{1})+1
\]
returns to 1 another also returns to 1. So when one of two numbers $n_{1}$ and $3^{r+1}2n_{i_{1}}(m_{1},t_{1})+1$ returns to 1 another also returns to 1. Let $m$ and $t$ be determined by
\[
n_{i}(m,t)=2n_{i}(m-1,t)+1\texttt{ and }n_{i}(m-1,t)=3^{r+1}n_{i_{1}}(m_{1},t_{1})
\]
where $0\leq i<m$ and one of two numbers $n_{i}(m-1,t)$ and $n_{i_{1}}(m_{1},t_{1})$ returns to 1. Then by Lemma~\ref{le:odd_seqs_n_m_m-1_2r+1}, when one of two numbers $n_{i}(m-1,t)$ and $n_{i_{1}}(m_{1},t_{1})$ returns to 1 another also returns to 1, and by Lemma~\ref{le:odd_seqs_n_m_m-1}, when one of the two Collatz's Sequences $S_{c}(m,t)$ and $S_{c}(m-1,t)$ returns to 1 another also returns to 1, so that when one of two numbers $n$ and $n_{i_{1}}(m_{1},t_{1})$ returns to 1 another also returns to 1.

For example, when $n_{i_{1}}(m_{1},t_{1})$ returns to 1, by Lemma~\ref{le:odd_seqs_n_m_m-1_2r+1} $S_{c}(m-1,t)$ returns to 1, and by Lemma~\ref{le:odd_seqs_n_m_m-1} $S_{c}(m,t)$ returns to 1, i.e., the odd number $3^{r+1}2n_{i_{1}}(m_{1},t_{1})+1$ returns to 1, then by Lemma~\ref{le:odd_seqs_3_4r_2n} the odd number $n_{1}=3\times4^{r}2n_{i_{1}}(m_{1},t_{1})+1$ returns to 1, so that the odd number $n$ returns to 1.

This completes the proof of the lemma.
\end{proof}

\begin{remark}\label{re:odd_seqs_n_m_m-1_r1}
By Lemma~\ref{le:odd_seq_r1} one half of positive integers $r_{1}$ satisfy $r_{1}=1$ and one quarter of positive integers $r_{1}$ satisfy $r_{1}=2$. Let the notation $P_{r}$ denote the process of mutual references or nested recursions to lemmas~\ref{le:odd_seqs_n_m_m-1}-\ref{le:odd_seqs_n_m_m-1_r1}.

In the proof of Lemma~\ref{le:odd_seqs_n_m_m-1}, one half of positive integers $r_{1}$ satisfy $r_{1}=1$ and directly make Lemma~\ref{le:odd_seqs_n_m_m-1} true, another half of positive integers $r_{1}$ continue the process $P_{r}$.

In the proof of Lemma~\ref{le:odd_seqs_n_m_m-1_r1}, one quarter of positive integers $r_{1}$ satisfy $r_{1}=2$ and directly make Lemma~\ref{le:odd_seqs_n_m_m-1_r1} true, and one half of positive integers $r_{1}$ satisfy $r_{1}=1$ and make the reference to Lemma~\ref{le:odd_seqs_n_m_m-1}, which is equivalent to that among these positive integers $r_{1}=1$ one half of them directly make Lemma~\ref{le:odd_seqs_n_m_m-1_r1} true and another half of them continue the process $P_{r}$. So it is equivalent to that one half of positive integers $r_{1}$ satisfy $r_{1}=1,2$ and directly make Lemma~\ref{le:odd_seqs_n_m_m-1_r1} true and another half of positive integers $r_{1}$ satisfy $r_{1}\neq2$ and continue the process $P_{r}$.

Thus Lemma~\ref{le:odd_seqs_n_m_m-1_r1} is equivalent to Lemma~\ref{le:odd_seqs_n_m_m-1} on the process $P_{r}$.
\end{remark}

\begin{lemma}\label{le:odd_seqs_3_4r_n}
When one of two numbers
\[
3\times4^{r}n_{i_{1}}(m_{1},t_{1})+1\texttt{ and }3^{r+1}n_{i_{1}}(m_{1},t_{1})+1
\]
returns to 1 another also returns to 1 where $r>0$ and $0\leq i_{1}\leq m_{1}$.
\end{lemma}

\begin{proof}
Let $n_{m_{0}^{*}}(m_{0}^{*},t_{0}^{*})=3\times4^{r}n_{i_{1}}(m_{1},t_{1})+1=4k+1$. Then there are
\[
(3n_{m_{i-1}^{*}}(m_{i-1}^{*},t_{i-1}^{*})+1)/2=2(3^{i+1}4^{r-i}n_{i_{1}}(m_{1},t_{1})+1)
\]
\[
\to n_{m_{i}^{*}}(m_{i}^{*},t_{i}^{*})=3^{i+1}4^{r-i}n_{i_{1}}(m_{1},t_{1})+1=4k_{i}+1\texttt{ for }i=1,2,\ldots,r-1,
\]
and
\[
n_{m_{r-1}^{*}+1}(m_{r-1}^{*},t_{r-1}^{*})=(3n_{m_{r-1}^{*}}(m_{r-1}^{*},t_{r-1}^{*})+1)/2=2(3^{r+1}n_{i_{1}}(m_{1},t_{1})+1).
\]

So when one of two numbers
\[
3\times4^{r}n_{i_{1}}(m_{1},t_{1})+1\texttt{ and }3^{r+1}n_{i_{1}}(m_{1},t_{1})+1
\]
returns to 1 another also returns to 1.

This completes the proof of the lemma.
\end{proof}

\begin{lemma}\label{le:odd_seqs_3_4r_2n}
When one of two numbers
\[
3\times4^{r}2n_{i_{1}}(m_{1},t_{1})+1\texttt{ and }3^{r+1}2n_{i_{1}}(m_{1},t_{1})+1
\]
returns to 1 another also returns to 1 where $r>0$ and $0\leq i_{1}\leq m_{1}$.
\end{lemma}

\begin{proof}
Let $n_{m_{0}^{*}}(m_{0}^{*},t_{0}^{*})=3\times4^{r}2n_{i_{1}}(m_{1},t_{1})+1=4k+1$. Then there are
\[
(3n_{m_{i-1}^{*}}(m_{i-1}^{*},t_{i-1}^{*})+1)/2=2(3^{i+1}4^{r-i}2n_{i_{1}}(m_{1},t_{1})+1)
\]
\[
\to n_{m_{i}^{*}}(m_{i}^{*},t_{i}^{*})=3^{i+1}4^{r-i}2n_{i_{1}}(m_{1},t_{1})+1=4k_{i}+1\texttt{ for }i=1,2,\ldots,r-1,
\]
and
\[
(3n_{m_{r-1}^{*}}(m_{r-1}^{*},t_{r-1}^{*})+1)/2=2(3^{r+1}2n_{i_{1}}(m_{1},t_{1})+1)
\]
\[
\to n_{i}(m_{r}^{*},t_{r}^{*})=3^{r+1}2n_{i_{1}}(m_{1},t_{1})+1\texttt{ where }0\leq i<m_{r}^{*}.
\]

So when one of two numbers
\[
3\times4^{r}2n_{i_{1}}(m_{1},t_{1})+1\texttt{ and }3^{r+1}2n_{i_{1}}(m_{1},t_{1})+1
\]
returns to 1 another also returns to 1.

This completes the proof of the lemma.
\end{proof}

\bibliographystyle{amsplain}

\end{document}